\documentclass[a4paper,10pt]{article}
\usepackage[utf8]{inputenc}

\usepackage{cmap}
\usepackage[T1]{fontenc}	

\usepackage[style=alphabetic,backend=bibtex]{biblatex}
\usepackage{microtype}

\usepackage[pdfauthor={Myriam Finster},
	colorlinks=true,
	linkcolor=black,
	citecolor=black,
	urlcolor=black,
	filecolor=black,
	pdftitle={Congruence Veech Groups},
	pdftex]{hyperref}
\usepackage{url}

\usepackage{tikz} 
\usetikzlibrary{matrix,arrows,calc,decorations.pathreplacing}

\pgfkeys{tikz/mymatrixenv/.style={decoration=brace,every left delimiter/.style={xshift=3pt},every right delimiter/.style={xshift=-3pt}}}
\pgfkeys{tikz/mymatrix/.style={matrix of math nodes,matrix anchor=center,left delimiter=(,right delimiter={)},inner sep=2pt,column sep=0.5em,row sep=0.5em,nodes={inner sep=0pt},text height=8pt}}
\pgfkeys{tikz/mymatrixbrace/.style={decorate,thick}}
\newcommand\mymatrixbraceoffseth{0.7em}
\newcommand\mymatrixbraceoffsetv{0.2em}

\newcommand*\mymatrixbraceright[4][m]{
    \draw[mymatrixbrace] ($(#1.north west)!(#1-#3-1.south west)!(#1.south west)-(\mymatrixbraceoffseth,0)$)
        -- node[left=2pt] {#4} 
        ($(#1.north west)!(#1-#2-1.north west)!(#1.south west)-(\mymatrixbraceoffseth,0)$);
}
\newcommand*\mymatrixbraceleft[4][m]{
    \draw[mymatrixbrace] ($(#1.north east)!(#1-#2-1.north east)!(#1.south east)+(\mymatrixbraceoffseth,0)$)
        -- node[right=2pt] {#4} 
        ($(#1.north east)!(#1-#3-1.south east)!(#1.south east)+(\mymatrixbraceoffseth,0)$);
}
\newcommand*\mymatrixbracetop[4][m]{
    \draw[mymatrixbrace] ($(#1.north west)!(#1-1-#2.north west)!(#1.north east)+(0,\mymatrixbraceoffsetv)$)
        -- node[above=2pt] {#4} 
        ($(#1.north west)!(#1-1-#3.north east)!(#1.north east)+(0,\mymatrixbraceoffsetv)$);
}
\newcommand*\mymatrixbracebottom[4][m]{
    \draw[mymatrixbrace] ($(#1.south west)!(#1-1-#3.south east)!(#1.south east)-(0,\mymatrixbraceoffsetv)$)
        -- node[below=2pt] {#4} 
        ($(#1.south west)!(#1-1-#2.south west)!(#1.south east)-(0,\mymatrixbraceoffsetv)$);
}

\usepackage{amssymb}
\usepackage{amsthm}
\usepackage{amsmath}
\usepackage{mathtools}
\usepackage{mathdots}
\usepackage{nicefrac}
\usepackage{enumitem}

\title{Congruence Veech Groups}
\author{Myriam Finster \\ \\ \normalsize Karlsruhe Institute of Technology (KIT) \\ \normalsize e-mail: myriam.finster@gmail.com}

\allowdisplaybreaks[2]

\newtheorem{prop}{Proposition}[section]
\newtheorem{thm}{Theorem}
\newtheorem{lemma}[prop]{Lemma}
\newtheorem{cor}[prop]{Corollary}

\theoremstyle{definition}
\newtheorem{defn}[prop]{Definition}

\newtheorem{rem}[prop]{Remark}
\newtheorem{ex}[prop]{Example}

\def\Aut{\mathrm{Aut}}
\def\Out{\mathrm{Out}}

\def\id{\mathrm{id}}

\def\SL{\mathrm{SL}}
\def\PSL{\mathrm{PSL}}

\def\GL{\mathrm{GL}}

\def\Stab{\mathrm{Stab}}

\def\Gal{\mathrm{Gal}}

\def\Aff{\mathrm{Aff}^+}
\def\Trans{\mathrm{Trans}}
\def\MCG{\mathrm{MCG}}

\def\pr{\mathrm{pr}}

\def\Z{\mathbb{Z}}
\def\N{\mathbb{N}}
\def\R{\mathbb{R}}
\def\H{\mathbb{H}}
\def\C{\mathbb{C}}
\def\Q{\mathbb{Q}}

\def\nl{\langle\!\langle\,}
\def\nr{\,\rangle\!\rangle}

\DeclareMathOperator{\wid}{width}
\DeclareMathOperator{\level}{level}

\DeclareMathOperator{\ab}{ab}
\DeclareMathOperator{\lift}{lift}
\DeclareMathOperator{\aff}{aff}
\DeclareMathOperator{\der}{der}
\DeclareMathOperator{\Kern}{ker}

\DeclareMathOperator{\lcm}{lcm}

\newcommand{\nom}[1]{\textit{#1}}

\begin{document}

\maketitle

\begin{abstract}
We study Veech groups of covering surfaces of primitive translation surfaces.
Therefore we define congruence subgroups in Veech groups of primitive translation surfaces using their action on the homology with entries in $\Z/a\Z$.
We introduce a congruence level definition and a property of a primitive translation surface which we call property $(\star)$. It guarantees that partition stabilising congruence subgroups of this level occur as Veech group of a translation covering.  

Each primitive surface with exactly one singular point has property $(\star)$ in every level. We additionally show that the surface glued from a regular $2n$-gon with odd $n$ has property $(\star)$ in level $a$ iff $a$ and $n$ are coprime. 
For the primitive translation surface glued from two regular $n$-gons, where $n$ is an odd number, we introduce a generalised Wohlfahrt level of subgroups in its Veech group. We determine the relationship between this Wohlfahrt level and the congruence level of a congruence group.
\end{abstract}

\section{Introduction}
A translation surface $\bar X$ is obtained by taking finitely many plane polygons and gluing their edges by translations in a way that leads to a connected, oriented surface without boundary.
The surface then carries a natural translation structure with finitely many cone angle singularities arising from the vertices of the polygons.
The derivatives of self-maps that are locally affine on the surface $X$ without the singularities form the Veech group $\Gamma(X)$. 
In 1989 Veech connected the properties of the geodesic flow on a translation surface to its Veech group (see \cite{Vee89}). 

An alternative way to define a translation surface is via a compact Riemann surface and a holomorphic one-form. 
This definition emphasises the connection between translation surfaces and the moduli space of Riemann surfaces. See e.g.\ \cite{Zo06} for an introduction to flat surfaces. 
In \cite{ER12} Ellenberg and McReynolds prove that many finite-index subgroups of $\SL_2(\Z)$ are Veech groups and that every algebraic curve $X/\bar\Q$ is birational over $\C$ to a Teichmüller curve.
  
Gutkin and Judge prove in \cite{GJ00} that if a translation surface $\bar Y$ covers a translation surface $\bar X$ via a map that is locally a translation, then their Veech groups $\Gamma(Y)$ and $\Gamma(X)$ are commensurate, i.e.\ $\Gamma(X) \cap \Gamma(Y)$ has finite index in $\Gamma(X)$ and in $\Gamma(Y)$. 
It is a consequence of \cite{Moe06} that if the translation surface $\bar X$ is primitive then $\Gamma(Y)$ is a finite index subgroup of $\Gamma(X)$. 
We address the question, which subgroups of the Veech group of a primitive translation surface $\bar X$ can be realised as Veech group of a covering surface $\bar Y$. 

Our starting point is the work of Gabriela Weitze-Schmithüsen in \cite{Sch05}.
There she proves that many congruence subgroups of $\SL_2(\Z)$ can be realised as Veech groups of Origamis, i.e.\ as Veech groups of translation coverings of the once-punctured torus.
We extend the definition of congruence groups to subgroups in the Veech group $\Gamma(X)$ for other primitive translation surfaces $\bar X$ and
generalise her results. 

We call a subgroup of $\Gamma(X)$ a congruence group of level $a$ if it contains each Veech group element that acts trivially on the homology with entries in $\Z/a\Z$.
We use a characteristic covering $\bar Y_a$ of $\bar X$ whose deck transformation group is the principal congruence group of level $a$ in $\Gamma(X)$ to define a property called property~$(\star)$ (see Definition \ref{def:star}). 
This Property~$(\star)$ guarantees the existence of affine maps stabilising a suitable subset of the singularities in $\bar Y_a$. It is the precondition to our first theorem which will be proven in Section~\ref{sec:main}.
\begin{thm}
 Let $a \geq 2$ and suppose that $\bar X$ has property $(\star)$ in level $a$.
 Furthermore, let $B = \{b_1, \dots, b_p\}$ be a partition of $H_1(\bar X, \Z/a\Z)$ and let 
 \begin{align*}
 \Gamma_B \coloneqq \{ \der(f) \mid & \, f \in \Aff(\bar X) \textrm{ and the map on } H_1(\bar X, \Z/a\Z) \\
 &\textrm{ induced by } f \textrm{ stabilises } B\}\,.
 \end{align*}
 Then there exists a translation covering $p \colon \bar Z \to \bar X$ with $\Gamma(Z) = \Gamma_B$.
\end{thm}
Every primitive translation surface whose pure Veech group equals its Veech group has property $(\star)$ in every level. 
This applies in particular to many of the surfaces introduced in Veech's article \cite{Vee89}. The surface glued from two regular $n$-gons, $\bar X_n$, for odd $n\geq 5$ and the surface glued from a regular $2n$-gon, $\bar X_{2n}$, for even $n \geq 4$ have exactly one singular point so their Veech groups and pure Veech groups (cf. Chapter \ref{sec:basics}) coincide.
The following theorem, proven in Section~\ref{sec:X2nStar}, gives examples of primitive surfaces with two singularities and property $(\star)$ in many but not in all levels.
\begin{thm}
 Let $n \geq 5$ be an odd number and $a \geq 2$.
 Then the translation surface $\bar X_{2n}$ has property $(\star)$ in level $a$ if and only if $\gcd(a,n)=1$.
\end{thm}

The group $\Gamma_B$ in Theorem~1 is a congruence group of level $a$. For congruence groups $\Gamma$ in $\SL_2(\Z)$ the Wohlfahrt level (see \cite{Wo63}) gives a complete characterisation of its congruence levels in terms of parabolic elements. 
In Section~\ref{sec:Wohlfahrt} we define a generalised Wohlfahrt level for finite index subgroups in the Veech group of $\bar X_n$ for odd $n\geq 5$. 
Here the congruence level and the Wohlfahrt level do not necessarily coincide but they still have a very strong connection:
\begin{thm}
 Let $\Gamma \leq \Gamma(X_n)$ be a congruence group, $b$ a minimal congruence level of $\Gamma$ and $a=\level(\Gamma)$ its generalised Wohlfahrt level. 
 Then $a \mid b$ and all prime numbers $p$ dividing $b$ also divide $a$. 
 However, $b$ does not have to divide $a$. 
\end{thm}
Finally we give a sufficient condition for a subgroup $\Gamma \leq \Gamma(X_n)$ to be not a congruence group and a explicit example of a non-congruence group in $\Gamma(X_5)$.

Most of the results stem from the author's PhD thesis \cite{Fin13}. 

\textbf{Acknowledgements}
I would like to thank the adviser of my thesis, Gabriela Weitze-Schmithüsen for many helpful discussions, inspiring suggestions and a lot of proofreading. I also thank the referee for his or her helpful comments. 

This work was partially supported by the Landesstiftung Baden-Württemberg within the program Juniorprofessuren-Programm.
\setcounter{thm}{0}

\section{Definitions and preliminaries}
\label{sec:basics}

In this section we shortly review the basic definitions used to state and prove our theorems. For a more detailed introduction to Veech groups see e.g.\! \cite{HS01} or \cite{Vor96}.

A \nom{translation surface} $\bar{X}$ is a connected, compact, $2$-dimensional, real manifold with a finite, non\-empty set $\Sigma(\bar X)$ of \nom{singular points} or \nom{singularities} together with a maximal $2$-dimensional atlas $\omega$ on $X = \bar{X} \setminus \Sigma(\bar X)$ such that all transition maps between the charts are translations. 
Furthermore, every singular point $s$ has an open neighbourhood $U$, not containing other singular points, such that there exists a continuous map $f_s \colon \hat U \to \hat V$ from $\hat U \coloneqq U \setminus \{s\}$ 
to a punctured open set $\hat V \subset \R^2$ that is compatible with $\omega$, i.e.\ $f_s \circ \varphi^{-1}$ is a translation for every $(U',\varphi) \in \omega$ with $U' \cap U \neq \emptyset$.

An alternative way to define a translation surface is by gluing finitely many Euclidean polygons via identification of edge pairs by translations. 
A translation structure on a torus for example can be obtained by gluing the parallel edges of a parallelogram.  
In this special case, no non-removable singularity arises from the vertices of the polygon and one has to add a removable singularity to meet the condition $\Sigma(\bar X) \neq \emptyset$. 
Especially in the situation where we glue the translation surface from a polygon the translation structure on $X$ is obvious, so we usually omit $\omega$ in the notation.

Let $\bar Y$ be a second translation surfaces.
We call a continuous map $p \colon \bar{Y} \to \bar{X}$ a \nom{translation covering} if $p^{-1}(\Sigma(\bar X)) = \Sigma(\bar Y)$ and $p|_Y \colon Y \to X$ is locally a translation. 
Furthermore, we call $\bar X$ the \nom{base surface} and $\bar Y$ the \nom{covering surface} of $p$. 
As translation surfaces are compact, a translation covering is a finite covering map in the topological sense, ramified at most over the singularities $\Sigma(\bar X)$. 
A translation surface $\bar X$ that does not admit a translation covering $\bar X \to \bar Y$ of degree $d > 1$ is called \nom{primitive}. 

An \nom{affine map} of a translation surface $\bar X$ is an orientation preserving homeomorphism $f$ on $\bar X$ with $f(\Sigma(\bar X)) = \Sigma(\bar X)$ that is affine on $X$, i.e.\ that can be written as $z \mapsto A z + b$ with $A \in \SL_2(\R)$ and $b \in \R^2$ in local coordinates. 
The translation vector $b$ depends on the local coordinates, whereas the \nom{derivative} $A$ is globally defined: the transition map between two charts is a translation and therefore does not change $A$. 

The \nom{affine group} $\Aff(\bar X)$ of $\bar X$ is the group of all affine maps on $\bar X$. 
The derivatives of the affine maps on $\bar X$ form the \nom{Veech group} $\Gamma(X) \leq \SL_2(\R)$ of $\bar X$. 
The \nom{projective Veech group} of $\bar X$ is the image of $\Gamma(X)$ in $\PSL_2(\R)$.
The affine maps with trivial derivative form the group of translations $\Trans(\bar X) \subseteq \Aff(\bar X)$.
An affine map is called \nom{pure} if it preserves the singularities pointwise (and not only setwise as usual).
The derivatives of all pure affine maps form the \nom{pure Veech group}. 

 As the set of singular points of a translation surface $\bar X$ is by definition nonempty, the fundamental group of $X$ is free of rank $n = 2g +(\nu -1)$ where $g$ is the genus of the surface and $\nu$ is the number of singular points. 
 We fix an isomorphism $\pi_1(X) \stackrel{\sim}{\longrightarrow} F_n$. 
 
 The affine group of $\bar X$ defines an outer action on $\pi_1(X)$.
 We define $\Aut_X(F_n)$ as the group of all automorphisms whose equivalence class in $\Out(F_n) \cong \Out(\pi_1(X))$ lies in the image of $\iota \colon \Aff(\bar X) \to \Out(\pi_1(X))$.
 Recall that the set of singular points $\Sigma(\bar X)$ is by definition nonempty. 
 By Lemma~5.2 in \cite{EG97} the affine group injects into the mapping class group $\MCG(\bar X)$ of $\bar X$ (which is the group of homeomorphisms on the surface $\bar X$ up to homotopy).
 It is a consequence of this lemma and the Dehn-Nielsen-Baer theorem (see e.g.\ Theorem~8.1 in \cite{FaMar12}) that the map $\iota$ is injective.
 Thus it induces a well-defined map 
 $$\Aut_X(F_n) \twoheadrightarrow \Out_X(F_n) \stackrel{\iota^{-1}}{\longrightarrow} \Aff(\bar X) \stackrel{\der}{\longrightarrow} \Gamma(X)$$
 that sends every $\gamma \in \Aut_X(F_n)$ to the derivative of the corresponding affine map. 
 The map is called $$\vartheta \colon \Aut_X(F_n) \to \Gamma(X)\,.$$
 For an affine map $f$ we call every preimage of $\iota(f)$ in $\Aut_X(F_n)$ a \nom{lift} of $f$ to $\Aut(F_n)$.
 
 In \cite{Moe06}, a translation surface is called primitive if it is not the covering of a translation surface of smaller genus. 
 By the Riemann-Hurwitz formula, the genus of the covering surface of a translation covering of degree $>1$ is always greater than the genus of the base surface, whenever the genus of the base surface is $> 1$. Thus the two definitions are equivalent for all surfaces of genus $>1$. 
 With our definition, there are non-primitive translation surfaces of genus $1$, whereas this is not the case in the alternative definition. 
 
 The following proposition states an important connection between the Veech group of the primitive base surface and the Veech group of the covering surface in a translation covering with primitive base surface. 
 It was proven for base surfaces of genus $>1$ in \cite{Moe06} Theorem~2.6 and for genus $1$ in \cite{Sch05} Proposition~3.3.

 \begin{prop}[see \cite{Moe06}, \cite{Sch05}]
 \label{prop:descend}
 Let $p \colon \bar Y \to \bar X$ be a translation covering with a primitive base surface.
 Every affine map on $\bar Y$ descends to $\bar X$.
\end{prop}

The relation between the two Veech groups can be specified even more precisely as follows:
\begin{prop}[see Korollar~6.22 in \cite{Fr08}]
\label{prop:Veech group by stabiliser}
 The Veech group element $A \in \Gamma(X)$ is contained in $\Gamma(Y)$ if and only if there exists a lift $\gamma$ of $A$ to $\Aut(\pi_1(X))$ such that $\gamma$ stabilises $\pi_1(Y) \leq \pi_1(X)$, i.e.\ 
 $$\Gamma(Y) = \vartheta (\Stab_{\Aut(\pi_1(X))}(\pi_1(Y))) \,.$$
\end{prop}

The proposition is a generalisation of Theorem~1 in \cite{Sch05}, using \cite{Sch08}.

 As stated above, every affine map $f$ on $\bar Y$ defines an automorphism $\gamma \in \Aut(\pi_1(Y))$, well-defined up to an inner automorphism of $\pi_1(Y) \leq \pi_1(X) \cong F_n$. 
 Proposition \ref{prop:descend} tells us that $f$ descends to $X$ and thus also defines a $\gamma' \in \Aut_X(F_n) \subseteq \Aut(F_n)$. 
 We define $\Aut_Y(F_n) \subseteq \Aut_X(F_n)$ as the set of those automorphisms.

\section{Translation coverings with congruence Veech groups}
\label{sec:main}
In the following, $\bar{X}$ is a primitive translation surface with $\nu \geq 1$ singularities and genus $g \geq 1$.
The fundamental group $\pi_1(\bar X)$ of the compact surface $\bar X$ is generated by $a_1, \dots, a_g, b_1, \dots, b_g$ where the $a_i, b_i$ belong to the $i$-th handle and fulfil the relation $a_1b_1a_1^{-1}b_1^{-1} \cdots a_gb_ga_g^{-1}b_g^{-1} = 1$. 
This relation does not hold in $\pi_1(X)$ because it describes a nontrivial path around the singularities. We amend the $a_i, b_i$ by paths $c_1, \dots, c_{\nu-1}$ to yield a basis of $\pi_1(X)$, where $c_i$ is a simple closed path around the $i$-th singularity of $X$. The group $\pi_1(X)$ is free of rank $n = 2 g + \nu -1$. We use the basis $\{a_1, b_1, \dots, a_g, b_g, c_1 , \dots, c_{\nu-1} \}$ to fix an isomorphism $\pi_1(X) \cong F_n$.

\subsection{Action on homology}
\label{sec:action on homology}
In $\SL_2(\Z)$, a congruence group of level $a$ is a subgroup of $\SL_2(\Z)$ that contains the kernel of the map $\bar \varphi_a \colon \SL_2(\Z) \to \SL_2(\Z/a\Z)$, obtained by sending each matrix entry to its residue modulo $a$. 
We know that the Veech group of the once-punctured torus $\bar E$ is $\SL_2(\Z)$. Furthermore, $H_1(\bar E, \Z/a\Z) \cong (\Z/a\Z)^2$ and the Veech group of $\bar E$ acts as $\SL_2(\Z/a\Z)$ on $H_1(\bar E, \Z/a\Z)$.
Thus the principal congruence group of level $a$ can equivalently be defined as the group of Veech group elements that act trivially on $H_1(\bar E, \Z/a\Z)$. 
We generalise this definition in a  straight forward way.

On a primitive translation surface $\bar X$, the affine group and the Veech group are isomorphic, as the translation group of $\bar X$ is trivial. 
We use this isomorphism and identify $\Gamma(X)$ with $\Aff(\bar X)$.
Thus $\Gamma(X)$ acts on the homology of $\bar X$ and we define the \nom{principal congruence group} of level $a$, $\Gamma(a)$, as the group of Veech group elements that act trivially on $H_1(\bar X, \Z/a\Z)$. 
As before, a \nom{congruence group} of level $a$ in $\Gamma(X)$ is a subgroup that contains $\Gamma(a)$.

\begin{center}
\begin{tikzpicture}
   \matrix (m) [matrix of math nodes, row sep=2.5em, column sep=2em, text height=2ex, text depth=0.25ex]
	{ & \pi_1(\bar X)\\
	F_n \cong \pi_1(X) & & H_1(\bar X, \Z) \cong \Z^{2g}\\
	& &  H_1(\bar X, \Z/a\Z) \cong (\Z/a\Z)^{2g}\\};
  \path[->, font=\footnotesize]
	(m-2-1) edge (m-1-2)
			edge node[auto] {$\ab$} (m-2-3)
			edge node[auto,swap] {$m_a$} (m-3-3)
	(m-1-2) edge node[auto] {$/[\pi_1(\bar X), \pi_1(\bar X)] $} (m-2-3)
	(m-2-3) edge node[auto] {$\pr_a$} (m-3-3);
\end{tikzpicture}
\end{center}

The action of $\Aff(\bar X)$ on the absolute homology $H_1(\bar X, \Z/a\Z) \cong (\Z/a\Z)^{2g}$ with entries in $\Z/a\Z$, can be derived from the outer action of the affine group on the fundamental group $\pi_1(X)$. 
We compose the group homomorphism $\pi_1(X) \to \pi_1(\bar X)$, given by $a_i \mapsto a_i$, $b_i \mapsto b_i$ and $c_i \mapsto 1$, with the abelianisation $\pi_1(\bar X) \to \pi_1(\bar X) / [\pi_1(\bar X), \pi_1(\bar X)] \cong H_1(\bar X, \Z)$ that maps the fundamental group of $\bar X$ to the absolute homology of $\bar X$ with integer coefficients.
The resulting homomorphism will be called $\ab \colon \pi_1(X) \to H_1(\bar X, \Z)$. 
The images of the $a_i$ and $b_i$ form a basis of $H_1(\bar X, \Z)$ and we use them to fix an isomorphism $H_1(\bar X, \Z) \cong \Z^{2g}$.
Next, we compose $\ab$ with the canonical projection $\pr_a \colon \Z^{2g} \to (\Z/a\Z)^{2g}$ and obtain the canonical homomorphism $m_a \colon \pi_1(X) \to H_1(\bar X , \Z/a\Z)$ from the fundamental group to the first homology of $\bar X$ with coefficients in $\Z/a\Z$.

Let $F_{2g} = \langle a_1, b_1, \dots, a_g, b_g\rangle \subset \pi_1(X)$ and let $F_{2g}^a$ be the set of all $a$-th powers of words in $F_{2g}$. 
Then a normal generating set for the kernel of $m_a$ is given by 
$$H \coloneqq \ker(m_a) = \nl \{c_1, \dots, c_{\nu-1}\} \cup [F_{2g},F_{2g}] \cup F_{2g}^a \nr \,.$$

\begin{lemma}
\label{lemma:resp_ma}
 Each element in $\Aut_X(F_n)$ respects $H = \Kern(m_a)$.
\end{lemma}
\begin{proof}
 Let $f \colon X \to X$ be an affine map and $\gamma \in \Aut_X(F_n)$ a lift of $f$. 
 The group $(\Z/a\Z)^{2g}$ is finite, so $H$ is of finite index in $F_n$ and it suffices to show that $\gamma$ maps the normal generators of $H$ to $H$.
  
 First, consider a generator $c_i$. As $f$ sends singular points to singular points, $\gamma(c_i) = w c_j w^{-1}$ for some $w \in F_n$ and $j \in \{1, \dots, \nu \}$. For $j \in\{1, \dots, \nu-1\}$ it is obvious that $\gamma(c_i) \in H$. 
 Furthermore, we see that 
 $$c_\nu = a_1b_1a_1^{-1}b_1^{-1} \cdots a_gb_ga_g^{-1}b_g^{-1} c_1^{-1} \cdots c_{\nu-1}^{-1} \in [F_{2g},F_{2g}] \cdot \langle c_1, \dots, c_{\nu-1}\rangle \subseteq H,$$ 
 thus $w c_\nu w^{-1} \in H$.
 
 For every $x,y \in F_{2g}$ we have that $m_a(\gamma(x y x^{-1}y^{-1})) = 0$, as $(\Z/a\Z)^{2g}$ is abelian, so $\gamma([F_{2g},F_{2g}]) \in H$. Furthermore, $a \cdot z = 0$ for every $z \in (\Z/a\Z)^{2g}$ implies that $\gamma(v^a) \in H$ for every $v \in F_{2g}$, so $\gamma(F_{2g}^a) \in H$.
\end{proof}

It follows from Lemma~\ref{lemma:resp_ma} that for every $\gamma \in \Aut_X(F_n)$ there is a unique homomorphism $\varphi_a(\gamma)$ that makes the diagram
\begin{center}
\begin{tikzpicture}
   \matrix (m) [matrix of math nodes, row sep=2.5em, column sep=4em, text height=1.5ex, text depth=0.25ex]
	{ F_n & F_n\\
	  (\Z/a\Z)^{2g} & (\Z/a\Z)^{2g}\\};
  \path[->, font=\footnotesize]
	(m-1-1) edge node[auto] {$\gamma$} (m-1-2)
			edge node[auto,swap] {$m_a$} (m-2-1)
	(m-1-2) edge node[auto] {$m_a$} (m-2-2)
	(m-2-1) edge node[auto,swap] {$\varphi_a(\gamma)$} (m-2-2);
\end{tikzpicture}
\end{center}
commutative. 
This defines an action $\varphi_a \colon \Aut_X(F_n) \to \Aut((\Z/a\Z)^{2g})$ of $\Aut_X(F_n)$ on $(\Z/a\Z)^{2g}$.
Observe that $\varphi_a$ actually defines an action $\bar{\varphi}_a \colon \Gamma(X) \to \Aut((\Z/a\Z)^{2g})$, since the homomorphism $\varphi_a(\gamma)$ does not depend on the chosen lift $\gamma$ of $f$ to $\Aut_X(F_n)$. 
The lift is unique up to an inner automorphism of $F_n$ and every inner automorphism clearly lies in the kernel of $\varphi_a$. 
Indeed, this is the standard action of $\Gamma(X) \cong \Aff(\bar X)$ on $H_1(\bar X, \Z/a\Z)$.

We derive the following characterisation of the principal congruence groups in $\Gamma(X)$. 
\begin{rem}
\label{rem:congruence group}
 The principal congruence group $\Gamma(a)$ of level $a$ in $\Gamma(X)$, where $\bar X$ is a primitive translation surface, is the kernel of the map $\bar{\varphi}_a$.
\end{rem}

\subsection{A characteristic covering}
\label{sec:charCovering}
If we use the map $m_a$ as monodromy map, we get a translation covering $\bar Y_a$ of $\bar X$ that is strongly related to the action of $\Gamma(X)$ on $H_1(\bar X, \Z/a\Z)$. 
For the definition of monodromy maps with our notation see \cite{Fin13} Section~1.2.
More details on monodromy maps can be found in \cite{Mir95}. 
So let $p_a\colon \bar Y_a \to \bar X$ be the translation covering of degree $a^{2g}$ defined by the monodromy map $m_a \colon F_n \twoheadrightarrow (\Z/a\Z)^{2g} \subseteq S_{a^{2g}}$, where $S_{a^{2g}}$ is the symmetric group on $a^{2g}$ elements. 
The covering map $p_a$ is unramified above $\Sigma(\bar X)$ because the $c_i$ lie in the kernel of $m_a$. 
Thus, by the Theorem of Riemann-Hurwitz, the genus of the surface $\bar Y_a$ is $g(\bar Y_a) = a^{2g} (g-1) +1$.

It is an immediate consequence of Lemma~\ref{lemma:resp_ma} that $p_a$ is a characteristic covering, i.e.\ all affine maps on $\bar X$ can be lifted to affine maps on $\bar Y_a$.

\begin{prop}
\label{prop:characteristic}
 The translation covering $p_a$ is characteristic and therefore $\Gamma(Y_a) = \Gamma(X)$.
\end{prop}
\begin{proof}
  Let $f \colon \bar X \to \bar X$ be an affine map. 
  By Proposition~\ref{prop:Veech group by stabiliser} we need to show that 
  a lift $\gamma$ of $f$ to $\Aut(F_n)$ preserves $\pi_1(Y_a)$ as subgroup of $\pi_1(X)$. 
  As $\pi_1(Y_a) = \ker(m_a) = H$, Lemma~\ref{lemma:resp_ma} tells us that all $\gamma \in \Aut_X(F_n)$ stabilise $\pi_1(Y_a)$.
\end{proof}

\begin{rem}
 Let $\mathcal{H}(d_1, \dots, d_\nu)$ be the stratum of translation surfaces containing $\bar X$. 
 Then Proposition~\ref{prop:characteristic} gives a characteristic translation covering whose covering surface lies in the stratum $\mathcal{H}(a^{2g}\, d_1, \dots, a^{2g}\,d_\nu)$ and has genus $g(\bar Y_a) = a^{2g} (g-1) +1$ for every $a \geq 2$.
\end{rem}

In the following we make the relation between the action of $\Gamma(X)$ on $H_1(\bar X, \Z/a\Z)$ and $\bar Y_a$ precise. Therefore we identify for every $s_i \in \Sigma(\bar X)$ the preimage set $p_a^{-1}(s_i) \subseteq \Sigma(\bar Y_a)$ with the elements in $H_1(\bar X, \Z/a\Z)$.

For a translation surface $\bar Z$ and a singularity $s \in \Sigma(\bar Z)$, we say that an element of $\pi_1(Z)$ is \nom{freely homotopic} to the singularity $s$, if it is nontrivial, freely homotopic to a simple path and can be freely homotoped into every neighbourhood of $s$. 

Let $\Sigma(\bar X) = \{s_1, \dots, s_\nu\}$. If we choose for each $i \in \{1,\dots, \nu\}$ a closed path $\hat c_i \in \pi_1(X)$ that is freely homotopic to $s_i$, we may construct a bijection $$\Sigma(\bar Y_a) \cong \Sigma(\bar X) \times H_1(\bar X, \Z/a\Z)$$ as follows. 
To simplify some arguments in Section~\ref{sec:ramCov} we require that $\hat c_i$ is a simple path in $\pi_1(Y_a)$. 
The $\hat c_i$ can be written as $\hat c_i = w_i c_i w_i^{-1}$ with $w_i \in \pi_1(X)$. 
As $p_a$ is unramified $\hat c_i$ lies in $\pi_1(Y_a) \subseteq \pi_1(X)$. 
Every simple closed path that is freely homotopic to $s \in \Sigma(\bar Y_a)$, can be written as $w' c_i w'^{-1}$ and consequently also as $w \hat c_i w^{-1}$ for suitable $i \in \{1, \dots, \nu\}$ and $w',w \in \pi_1(X)$. 
Two elements $w \hat c_i w^{-1}$ and $w' \hat c_j w'^{-1}$ in $\pi_1(Y_a)$ are homotopic to the same singularity iff $i = j$ and $m_a(w) = m_a(w')$. 
We use this to identify $s$ with $m_a(w)\in (\Z/a\Z)^{2g}$. 
For every $s_i \in \Sigma(\bar X)$ we obtain the bijection 
$$\begin{array}{rcl}
 \Sigma(\bar Y_a) \supseteq p_a^{-1}(s_i) &\stackrel{\sim}{\longrightarrow} & H_1(\bar X, \Z/a\Z)\\
 s &\mapsto & m_a(w)\;, \textrm{ with } w \hat c_i w^{-1} \textrm{ freely homotopic to } s\,.
\end{array}$$

Of course this identification depends on the initially chosen path $\hat c_i$. 
But $\pi_1(Y_a) = \ker(m_a)$ implies that if $\hat c_i$ and $\hat c_i'$ are homotopic to the same singularity of $\bar Y_a$, i.e. differ by conjugation with an element in $\pi_1(Y_a)$, then they define the same map. 
Altogether the bijections define a map 
$$\tilde m_a \colon \Sigma(\bar Y_a) \twoheadrightarrow H_1(\bar X, \Z/a\Z)$$
and more precisely a bijection
$$\begin{array}{rcll}
\Sigma(\bar Y_a) & \stackrel{\sim}{\longrightarrow} & \Sigma(\bar X) \times H_1(\bar X, \Z/a\Z) \\
s & \mapsto & (p_a(s), \tilde m_a(s)) &.
\end{array}$$
Thus every singularity of $\bar Y_a$ is uniquely defined by a pair $(i, z)$, where $i \in \{1, \dots, \nu\}$ and $z \in H_1(\bar X, \Z/a\Z) \cong (\Z/a\Z)^{2g}$. 

\begin{figure}[htbp]
 \centering
\begin{tikzpicture}

\tikzset{
    myarrow/.style={->, >=latex', shorten >=1pt, thick}
} 
\tikzset{
	mydarrow/.style={<->, >=latex', shorten >=1pt, thick}
}

 \draw (0,0) -- node[below] {$a$} (2,0); 
 \draw (2,0) -- node[below] {$b$} (4,0);
 
 \draw (0,4) -- node[above] {$a$} (2,4); 
 \draw (2,4) -- node[above] {$b$} (4,4);
 
 \draw (0,0) -- node[left] {$c$} (0,2);
 \draw (0,2) -- node[left] {$d$} (0,4);
 
 \draw (4,0) -- node[right] {$c$} (4,2);
 \draw (4,2) -- node[right] {$d$} (4,4);

 \draw[dashed] (0,2) -- (4,2);
 \draw[dashed] (2,0) -- (2,4);

 \draw (-4,1) 
  -- node[below] {$y^{-1}$} ++(2,0) 
  -- node[right] {$x$} ++(0,2) 
  -- node[above] {$y$} ++ (-2,0) 
  -- node[left] {$x^{-1}$} ++ (0,-2);

 \fill (-3,2) circle (2pt);
 \draw[myarrow] (-3,2)
	.. controls (-2.5,2) and (-2.5,2.5) .. (-2,2.5);
 \draw[myarrow] (-4,2.5)
	.. controls (-3.5,2.5) and (-3.5,2.5) .. (-3.5,3);
 \draw[myarrow] (-3.5,1)
	.. controls (-3.5,1.5) and (-3.5,1.5) .. (-4,1.5);
 \draw[myarrow] (-2,1.5)
	.. controls (-2.5,1.5) and (-2.5,1.5) .. (-2.5,1);
 \draw[myarrow] (-2.5,3) 
	.. controls (-2.5,2.5) and (-3,2.5) .. (-3,2);
	
 \fill (1,1) circle (2pt);
 \draw[mydarrow] (1,1) -- node[above right]{$w$} (3,1);
 
 \draw[myarrow] (3,1)
  .. controls (3.5,1) and (3.5,1.5) .. (4,1.5);
 \draw[myarrow] (0,1.5)
	.. controls (0.5,1.5) and (0.5,1.5) .. (0.5,2);
 \draw[myarrow] (0.5,2)
	.. controls (0.5,2.5) and (0.5,2.5) .. (0,2.5);
 \draw[myarrow] (4,2.5)
	.. controls (3.5,2.5) and (3.5,2.5) .. (3.5,2);
 \draw[myarrow] (3.5,2) 
	.. controls (3.5,1.5) and (3,1.5) .. (3,1);

\draw (1,1) node[below, font=\tiny] {$\begin{pmatrix}
                              0 \\ 0
                             \end{pmatrix}$
};
 \draw (3,1) node[below, font=\tiny] {$\begin{pmatrix}
                              1 \\ 0
                             \end{pmatrix}$
};
 \draw (1,3) node[font=\tiny] {$\begin{pmatrix}
                              0 \\ 1
                             \end{pmatrix}$
};
 \draw (3,3) node[font=\tiny] {$\begin{pmatrix}
                              1 \\ 1
                             \end{pmatrix}$
};
	
 \fill (3.9,1.9) rectangle (4.1,2.1);
 \draw (4,2) node[right, font=\tiny] {$\begin{pmatrix}
                              1 \\ 0
                             \end{pmatrix}$
};
\fill (-0.1,1.9) rectangle (0.1,2.1);
\draw (0,2) node[left, font=\tiny] {$\begin{pmatrix}
                              1 \\ 0
                             \end{pmatrix}$
};

\fill[black!60] (1.9,1.9) rectangle (2.1,2.1);
\draw (2,2) node[above left, font=\tiny] {$\begin{pmatrix}
                              0 \\ 0
                             \end{pmatrix}$
};

\fill[black!40] (1.9,3.9) rectangle (2.1,4.1);
\draw (2,4) node[above, font=\tiny] {$\begin{pmatrix}
                              0 \\ 1
                             \end{pmatrix}$
};
\fill[black!40] (1.9,-0.1) rectangle (2.1,0.1);
\draw (2,0) node[below, font=\tiny] {$\begin{pmatrix}
                              0 \\ 1
                             \end{pmatrix}$
};

\fill[black!20] (-0.1,-0.1) rectangle (0.1,0.1);
\draw (0,0) node[below, font=\tiny] {$\begin{pmatrix}
                              1 \\ 1
                             \end{pmatrix}$
};
\fill[black!20] (3.9,-0.1) rectangle (4.1,0.1);
\draw (4,0) node[below, font=\tiny] {$\begin{pmatrix}
                              1 \\ 1
                             \end{pmatrix}$
};
\fill[black!20] (-0.1,3.9) rectangle (0.1,4.1);
\draw (0,4) node[above, font=\tiny] {$\begin{pmatrix}
                              1 \\ 1
                             \end{pmatrix}$
};
\fill[black!20] (3.9,3.9) rectangle (4.1,4.1);
\draw (4,4) node[above, font=\tiny] {$\begin{pmatrix}
                              1 \\ 1
                             \end{pmatrix}$
};

\end{tikzpicture}
\caption{Singularities in $\bar Y_2$, identified with $H_1(\bar E, \Z/2\Z) \cong (\Z/2\Z)^2$.}
\label{fig:YaSing}
\end{figure}

\begin{ex}
 The simplest example of this identification is shown in Figure~\ref{fig:YaSing}. There, the once-punctured torus $E$ is used as primitive base surface, glued of a unit square with the identified vertices as unique singularity. The centre of the square is used as base point of the fundamental group and the horizontal closed path $x$ and the vertical closed path $y$ through the centre as free generating set. Then the surface $Y_2$ consists of four copies of $E$, labelled by 
 $\left(\begin{smallmatrix} 0 \\ 0 \end{smallmatrix}\right)$,
 $\left(\begin{smallmatrix} 1 \\ 0 \end{smallmatrix}\right)$,
 $\left(\begin{smallmatrix} 0 \\ 1 \end{smallmatrix}\right)$ and
 $\left(\begin{smallmatrix} 1 \\ 1 \end{smallmatrix}\right)$. 
 A simple closed path, homotopic to the singularity of $\bar E$, is given by $x y x^{-1} y^{-1}$. We use it as $\hat c_1$. 
 Furthermore, we choose the centre of the copy labelled by 
 $\left(\begin{smallmatrix} 0 \\ 0 \end{smallmatrix}\right)$ 
 as base point of $\pi_1(Y_2)$. 
 Figure~\ref{fig:YaSing} shows the resulting identification of the singularities in $Y_2$ with $(\Z/2\Z)^2$ through $\tilde m_2$. 
 In addition, a path $w \hat c_1 w ^{-1}$ with $m_2(w) = \left(\begin{smallmatrix} 1 \\ 0 \end{smallmatrix}\right)$ is drawn to demonstrate the correlation to the singularity labelled by $\left(\begin{smallmatrix} 1 \\ 0 \end{smallmatrix}\right)$.  
 \end{ex}

Now we describe the action of an affine map $f \in \Aff(\bar Y_a)$ on the singularities of $\bar Y_a$ via the identification with $\Sigma(\bar X) \times H_1(\bar X, \Z/a\Z)$.

\begin{prop}
\label{prop:ActionOnSingbyHom}
 Let $f$ be an affine map of $\bar Y_a$ with derivative $A$. Furthermore, let $\gamma$ be any lift of $f$ to $\Aut_{Y_a}(F_n)$. 
 As affine maps send singular points to singular points, there are $j_i \in \{1, \dots, \nu\}$ and $v_i \in F_n$ such that $\gamma(\hat c_i) = v_i \hat c_{j_i} v_i^ {-1}$. 

Then $f(i,z) = (j_i, \bar \varphi_a(A) \cdot z + z_i)$ where $z_i = m_a(v_i)$.
\end{prop}
\begin{proof}
By definition $f(i,0) = (j_i, z_i)$.

Let $w \in F_n$ with $m_a(w) = z$, i.e.\ such that $w \hat c_i w ^{-1}$ is freely homotopic to the singularity $(i,z)$.
Then 
$\gamma(w \hat c_i w ^{-1}) = \gamma(w) v_i \hat c_{j_i} v_i^{-1} \gamma(w)^{-1}$ and 
$$m_a(\gamma(w) v_i) = (\varphi_a(\gamma))(m_a(w)) + m_a(v_i) = \bar \varphi_a(A) \cdot z + z_i\,.\qedhere$$ 
\end{proof}

Now suppose that we manage to find a section $\aff \colon \Gamma(Y_a) \to \Aff(Y_a)$ of the map $\der$ and that we further accomplish to choose the $\hat c_i \in \pi_1(X)$ such that $z_i = 0$ for every $f$ in the image of $\aff$ and all $i \in \{1, \dots, \nu \}$. 
Then the action of $\Gamma(X) = \Gamma(Y_a)$ via $\aff$ on $\Sigma(\bar Y_a) \cong \Sigma(\bar X) \times H_1(\bar X, \Z/a\Z)$ is compatible with the projection of $\Sigma(\bar X) \times H_1(\bar X, \Z/a\Z)$ to $H_1(\bar X, \Z/a\Z)$. 
Proposition~\ref{prop:ActionOnSingbyHom} states that the induced action on $H_1(\bar X, \Z/a\Z)$ equals the standard action of $\Gamma(X)$ on $H_1(\bar X, \Z/a\Z)$.  
However, we do not need all $z_i$ to be equal to zero, if we want to recover the action of $\Gamma(X)$ on the homology elements in its action on the singularities of $Y_a$ via $\aff$. 
What we really need is summarised in the following definition:

\begin{defn}
\label{def:star}
 Let $\bar X$ be a primitive translation surface with singularities $\{s_1, \dots, s_\nu\}$ and $a \geq 2$. 
 Furthermore, let $\{A_j \mid j \in J\}$ be a system of generators of $\Gamma(X)$.  
 The surface is said to have property $(\star)$ in level $a$ iff 
  there exists a $\mu \in \{1, \dots, \nu\}$ and
  $S = \{\hat s_1, \dots, \hat s_\mu\} \subseteq \Sigma(\bar Y_a)$ such that $p_a (\hat s_i) = s_i$ for $i \in \{1, \dots, \mu\}$ and 
  such that for every $j \in J$ there is an affine map $f_j \in \Aff(\bar Y_a)$ with $\der(f_j) = A_j$ and $f_j(S) = S$.
\end{defn}

The condition $f(S) = S$ is closed under composition, so property $(\star)$ assures that for every $A \in \Gamma(Y_a) = \Gamma(X)$ there exists an $f \in \Aff(\bar Y_a)$ with $f(S) = S$.  
Now suppose that $f$, $f' \in \Aff(\bar Y_a)$ both fix $S$ and have the same derivative. Then $f' \circ f^{-1} \in \Trans(\bar Y_a)$. 
As $\bar X$ is primitive, every affine map on $\bar Y_a$ descends to $\bar X$ via $p_a$ and with $\Trans(\bar X) = \{\id\}$ we obtain $f' \circ f^{-1} \in \Gal(\bar Y_a / \bar X)$. 
As $p_a \colon \bar Y_a \to \bar X$ is unramified, $f' \circ f^{-1}$ is either without fixed points or equals the identity map. 
Let $f(\hat s_i) = \hat s_j$. As $f$ and $f'$ have the same derivative, they both descend to the same map on $X$, mapping $s_i$ to $s_j$. This implies $f'(\hat s_i) = \hat s_j$, so $f' \circ f^{-1}$ has the fixed point $\hat s_j$ and thereby $f = f'$.

Altogether this proves, that the set $\{f_j \mid j \in J\}$ defines via $A_j \mapsto f_j$ a section $\aff \colon \Gamma(Y_a) \to \Aff(\bar Y_a)$.

\begin{rem}
 Let $\bar X$ be a primitive translation surface with property $(\star)$ in level $a$. 
 We choose the $\hat c_i \in \pi_1(Y_a) \subseteq \pi_1(X)$ to be freely homotopic to $\hat s_i$ for $i \in \{ 1, \dots, \mu\}$.
 According to property $(\star)$ and Proposition \ref{prop:ActionOnSingbyHom} every element in the image of $\aff$ respects $\Sigma' \coloneqq p_a^{-1}(p_a(S)) \cong p_a(S) \times H_1(\bar X,\Z/a\Z)$. Thus the action of $\Gamma(Y_a)$ on $\Sigma(\bar Y_a)$ via $\aff$ can be restricted to an action on $\Sigma'$. 
 As $z_i = 0$ for every $i \in \{1, \dots, \mu\}$, this action induces via $\tilde{m}_a \colon \Sigma(\bar Y_a) \to H_1(\bar X, \Z/a\Z), (i,z) \mapsto z$ an action of $\Gamma(X) = \Gamma(Y_a)$ on $H_1(\bar X, \Z/a\Z)$ which equals the standard action $\bar \varphi_a$ of $\Gamma(X)$ on $H_1(\bar X, \Z/a\Z)$.
\end{rem}
 
\begin{rem}
 Primitive translation surfaces with exactly one singular point and more generally, translation surfaces whose pure Veech group equals its Veech group obviously have property $(\star)$ in every level: 
 use $\mu = 1$ and chose an arbitrary $\hat s_1 \in \Sigma(\bar Y_a)$.
 As the unramified covering $p_a \colon \bar Y_a \to \bar X$ is normal, $\Trans(\bar Y_a)$ acts transitively on the set $\Sigma' = p_a^{-1}(p_a(\hat s_1))$. Every $f \in \Aff(\bar X)$ fixes $p_a(\hat s_1)$, thus it has a lift fixing $\hat s_1$.
\end{rem}

\subsection{A ramified covering}
\label{sec:ramCov}

If $\bar X$ has property $(\star)$ in level $a$, then the derivatives of affine maps that fix the set $\{(i,z) \mid  i \in \{1, \dots, \mu \} \} \subseteq \Sigma(\bar Y_a)$ for every $z \in H_1(\bar X, \Z/a\Z)$ form the principal congruence group of level $a$. 
In this section we make use of that property. 
For every congruence group of level $a$ that equals a set $\Gamma_B$ as defined in Definition~\ref{def:GammaB} and thus is the preimage set of the stabiliser of a partition of $H_1(\bar X, \Z/a\Z)$, we give a translation covering of $\bar X$ with Veech group $\Gamma_B$. 
To achieve this we define a ramified covering of $\bar Y_a$ whose ramification behaviour is defined using the (via $\tilde m_a$) induced partition of $\Sigma'$.

\begin{defn}
\label{def:GammaB}
 Let $B = \{ b_1, \dots, b_p \}$ be a partition of $(\Z/a\Z)^{2g}$, i.e.\ $b_i \cap b_j = \emptyset$ for $i \neq j$ and $\bigcup_{1=1}^p \, b_i = (\Z/a\Z)^{2g}$.
Define
\begin{align*}
\Gamma_B &\coloneqq \{A \in \Gamma(X) \mid \forall i\in \{1, \dots, p\} \exists j\in \{1, \dots, p\} : \bar\varphi_a(A)(b_i) = b_j\}\,.
\end{align*}
\end{defn}

\begin{thm}
\label{thm:GammaBAsVeechGroup}
 Let $a \geq 2$ and suppose that $\bar X$ has property $(\star)$ in level $a$.
 Furthermore, let $B = \{b_1, \dots, b_p\}$ be a partition of $(\Z/a\Z)^{2g}$. Then there exists a translation covering $p \colon \bar Z \to \bar X$ with $\Gamma(Z) = \Gamma_B$.
\end{thm}

The proof of this theorem is the goal of this section. 
Thus, from now on we assume that $a \geq 2$ and that $\bar X$ is a primitive translation surface with property $(\star)$ in level $a$.

We choose $p$ different natural numbers 
$r_{1}, \dots, r_{p}$ 
greater than $1$ and define three subsets of $H = \pi_1(Y_a)$: 
the first one is the set $C_X$ of all simple closed curves that are freely homotopic to a singularity in $\bar X$:
$$C_X \coloneqq \{w \hat c_i w^{-1} \mid w \in F_n, i \in \{1, \dots, \nu\} \, \} \subseteq H \,.$$
As $p_a$ is unramified, $C_X$ is also the set of all closed curves, freely homotopic to a singularity in $\bar Y_a$.

The second set, $C_{B,\mu}$, contains closed curves that are a power of elements from $C_X$, chosen in accordance with the partition $B$: 
\begin{align*}
  C_{B,\mu} &\coloneqq \{ w \hat c_i^{r_\varsigma} w^{-1} \mid w \in F_n,  i \in \{1, \dots, \mu\} \textrm{ and } \varsigma \textrm{ is such that } m_a(w) \in b_\varsigma \,\} \\
 & \quad \cup \quad \{ w \hat c_i w^{-1} \mid w \in F_n , i \in \{\mu+1, \dots, \nu\} \, \}
 \qquad \qquad \qquad \qquad \quad \subseteq H \,.
\end{align*}

Following the proof of Lemma~6.5 in \cite{Sch05}, we now prove a characterisation of paths winding (several times) around a singularity and lying in $N_{B,\mu} \coloneqq \langle C_{B,\mu} \rangle$.

\begin{lemma}
 \label{lemma:inNBmu}
  Let $w \in F_n$ with $m_a(w) \in b_\varsigma$ and $i \in \{1, \dots, \mu\}$ then 
  $$w {\hat c_i}^l w^{-1} \in N_{B,\mu} = \langle C_{B,\mu} \rangle \Leftrightarrow r_\varsigma \mid l\,.$$
\end{lemma}
\begin{proof} 
 As  $w \hat c_i^{r_\varsigma} w^{-1} \in C_{B,\mu}$, $r_\varsigma \mid l$ obviously implies $w {\hat c_i}^l w^{-1} \in N_{B,\mu}$. 
 
 To prove the reverse implication, let $p_\infty \colon Y_\infty \to X$ be the unramified covering defined by the subgroup $N_X \coloneqq \langle C_X \rangle$, i.e.\ the normal unramified covering with monodromy map $m \colon F_n \to F_n / N_X$. 
 By lifting the charts from $X$ to $Y_\infty$, the surface $Y_\infty$ becomes an infinite translation surface. 
 We show in the following that $N_X = \pi_1(Y_\infty)$ is freely generated by the set $S$ defined in Equation~\ref{eqn:S}.
 
 For every $i \in \{1, \dots, \nu \}$, the path $c_i \in \pi_1(X)$ is freely homotopic to the singularity $s_i$ in $\bar X$. 
 Hence if we develop the path $c_i$ along appropriate charts in $\R^2$, then we get a closed curve that has winding number $\kappa$ around an innermost point iff the singularity $s_i$ has multiplicity $\kappa$. 
 The path $\hat c_i$ is contained in $N_X$, so $m(\hat c_i) = 1_{F_n/N_X}$, thus $m(c_i) = 1_{F_n/N_X}$. 
 This implies that the path $c_i$ also describes a closed path in $Y_\infty$ with finite winding number $\kappa$ when projected to $\R^2$. 
 The same is obviously true for every conjugate of $c_i$.
 Thus the metric completion $Y'_\infty$ of $Y_\infty$ adds a finite angle singularity to the translation structure of $Y_\infty$ for every singularity $s_i \in \Sigma(\bar X)$ and for every $k \in F_n/ N_X$. 
 The covering $Y_\infty \to X$ extends to an unramified covering map $Y'_\infty \to \bar X$ and we get a commutative diagram \\ 
 \vspace{-1cm}
\begin{center}
\begin{tikzpicture}
   \matrix (m) [matrix of math nodes, row sep=2.5em, column sep=2em, text height=1.5ex, text depth=0.25ex]
	{Y_\infty & Y'_\infty\\
	X & \bar X\\};
  \path[right hook->, font=\footnotesize]
	(m-1-1) edge node[auto] {$\beta$} (m-1-2)
	(m-2-1) edge node[auto] {$\alpha$} (m-2-2);
  \path[->, font=\footnotesize]
  	(m-1-1) edge (m-2-1)
	(m-1-2) edge (m-2-2);
\end{tikzpicture}
\end{center}
\vspace{-0.3cm}
with inclusions $\alpha \colon X \hookrightarrow \bar X$ and $\beta \colon Y_\infty \hookrightarrow Y'_\infty$.
The maps $\alpha$ and $\beta$ induce maps $\alpha_\ast \colon \pi_1(X) \to \pi_1(\bar X)$ and $\beta_\ast \colon \pi_1(Y_\infty) \to \pi_1(Y'_\infty)$ and a commutative diagram of the fundamental groups:\\ 
 \vspace{-0.5cm}
 \begin{center}
\begin{tikzpicture}
   \matrix (m) [matrix of math nodes, row sep=2.5em, column sep=2em, text height=1.5ex, text depth=0.25ex]
	{N_X = \pi_1(Y_\infty) & \pi_1(Y'_\infty)\\
	\pi_1(X) & \pi_1(\bar X)\\};
  \path[right hook->, font=\footnotesize]
	(m-1-1) edge (m-2-1)
	(m-1-2) edge (m-2-2);
  \path[->>, font=\footnotesize]
  (m-1-1) edge node[auto] {$\beta_\ast$} (m-1-2)
  (m-2-1) edge node[auto] {$\alpha_\ast$} (m-2-2);
\end{tikzpicture}
\end{center}
 The surfaces $X$ and $\bar X$ and thereby also $Y_\infty$ and $Y'_\infty$ differ only in a discrete set of points. 
 Hence $\alpha_\ast$ and $\beta_\ast$ are surjective. 
 Furthermore, $\alpha_\ast(c_i) = 1$ thus $\alpha_\ast(\hat c_i) = 1$ and $\alpha_\ast(N_X) = 1$. 
 This implies that $\pi_1(Y'_\infty) = \beta_\ast(N_X) = \alpha_\ast(N_X) = \{1\}$ is trivial. 
 Thus $Y'_\infty$ is the universal covering of $\bar X$ and in particular simply connected. 
 So the genus of $Y'_\infty$ is $0$. 
 Consequently, $Y_\infty$ is homeomorphic to a plane with a discrete subset of points removed. 
 Thus $\pi_1(Y_\infty)$ is freely generated by a set that contains a simple closed path around each singularity of $Y_\infty$. 
 For every $s \in \Sigma(\bar X)$ and every element in $F_n/N_X$ there is a singularity in $Y_\infty$. 

 As $N_X \subseteq H$, $Y_\infty$ is a covering of $Y_a$. 
 The elements $\hat c_i$ are by definition simple closed paths in $Y_a$ thus they are simple in $Y_\infty$. 
 If we choose a preimage $v_h$ of every $h \in F_n/N_X$, then the following set is a free generating set of $\pi_1(Y_\infty)$:
 \begin{align}
 \label{eqn:S}
  S \coloneqq \{v_h \hat c_i v_h^{-1} \mid h \in F_n/N_X, i \in \{1, \dots, \nu \} \, \} 
 \,.
 \end{align}
 For $i \in \{1, \dots, \mu\}$ and $h \in F_n/N_X$ we define the homomorphism
 $$
  \varphi_{i,h}\colon \pi_1(Y_\infty) \to (\Z,+), w \mapsto \sharp_{v_h \hat c_i v_h^{-1}} (w) \,.
 $$
by sending the free generator $v_h \hat c_i v_h^{-1}$ to $1$ and the remaining generators to $0$.

For $w \in F_n$ with $m(w) = k$ it follows that
\begin{align*}
 \varphi_{i,h}(w{\hat c_j}^lw^{-1}) &= \varphi_{i,h}(w v_k^{-1}) + l \cdot \varphi_{i,h}(v_k \hat c_j v_k^{-1}) - \varphi_{i,h}(w v_k^{-1})\\
	&= \left \{ \begin{array}{ll}
	            0&, \textrm{if } k \neq h \textrm{ or } j \neq i\\ 
	            l&, \textrm{if } k=h \textrm{ and } j = i
	           \end{array} \right. \,.
\end{align*}	           
As $N_X \subseteq H$ the map $m_a$ factors through $m$: 
  \begin{center}
\begin{tikzpicture}
   \matrix (m) [matrix of math nodes, row sep=1em, column sep=2em, text height=1.5ex, text depth=0.25ex]
	{ F_n & & F_n/H \cong (\Z/a\Z)^{2g}\\
	  & F_n/N_X\\};
  \path[->, font=\footnotesize]
	(m-1-1) edge node[auto] {$m_a$} (m-1-3)
			edge node[auto,swap] {$m$} (m-2-2)
	(m-2-2) edge node[auto,swap] {$\phi$} (m-1-3);
\end{tikzpicture}
\end{center}
As $C_{B,\mu} \subseteq N_X$ we have $N_{B,\mu} \subseteq N_X = \pi_1(Y_\infty)$. 
    
If $w, w' \in F_n$ with $m(w) = m(w') = k$, then $m_a(w') = \phi(m(w')) = \phi(m(w)) = m_a(w)$. Thus all $w \in m^{-1}(k)$ are mapped into the same $b_\varsigma \subseteq (\Z/a\Z)^{2g}$. 
We conclude for $i \in \{1, \dots, \mu\}$, $k \in F_n/N_X$ and $\varsigma$ such that $\phi(k) \in b_\varsigma$ that
$$ \varphi_{i,k}(C_{B,\mu}) = \varphi_{i,k}(\{w \hat c_i^{r_\varsigma} w^{-1} \mid w \in F_n, m(w) = k\}) = \{r_\varsigma\}\,,$$
implying
$\varphi_{i,k}(N_{B,\mu}) =  \langle r_\varsigma \rangle$. 

Now let $w {\hat c_i}^l w^{-1} \in N_{B,\mu}$ with $m(w) = k$ and let $\varsigma$ be such that $m_a(w) \in b_\varsigma$. 
Then $\varphi_{i,k}(w {\hat c_i}^l w^{-1}) = l \in \langle r_\varsigma \rangle$ thus $r_\varsigma \mid l$. 
\end{proof}

In analogy to Lemma 6.7 in \cite{Sch05} we prove that the stabiliser of $N_{B,\mu}$ in $\Aut_X(F_n)$ is the set 
\begin{align*}
    G_{B,\mu} = 
  \{\gamma \in \Aut_{X}(F_n) \mid & \;
  \forall i \in \{1, \dots, \mu\} \, \exists j_i \in \{1, \dots, \mu\} \textrm{ and } v_i \in F_n : \\
  & \; \gamma(\hat c_i) = v_i \hat c_{j_i} v_i^{-1},\, m_a(v_i) = 0 \textrm{ and }  \forall w \in F_n:\\
  & \; m_a(w) \textrm{ and } m_a(\gamma(w)) \textrm{ are in the same } b_\varsigma\}\,.
\end{align*}

\begin{lemma}
\label{lemma:stab=stab_mu}
 $G_{B,\mu} = \Stab_{\Aut_X(F_n)}(C_{B,\mu}) = \Stab_{\Aut_X(F_n)}(N_{B,\mu})$
\end{lemma}
\begin{proof}
 The inclusion $\Stab_{\Aut_X(F_n)}(C_{B,\mu}) \subseteq \Stab_{\Aut_X(F_n)}(N_{B,\mu})$ is trivial. 
 
 Now let $\gamma \in \Stab_{\Aut_X(F_n)}(N_{B,\mu})$.
 We start by showing the following claim: 
 if $i \in \{1, \dots, \mu \}$ then $\gamma(\hat c_i) = v_i \hat c_{j_i} v_i^{-1}$ where $j_i \in \{1, \dots, \mu\}$ (and $v_i \in F_n$). 
 Furthermore, if $i \in \{\mu+1, \dots, \nu\}$ then $\gamma(\hat c_i) = v_i \hat c_{j_i} v_i^{-1}$ with $j_i \in \{\mu+1, \dots, \nu\}$ (and $v_i \in F_n$).
 
 So at first let $i \in \{1, \dots, \mu\}$ and suppose that $\gamma(\hat c_i) = v_i \hat c_j v_i^{-1}$ where $j > \mu$. 
 Then $v_i \hat c_j v_i^{-1} \in N_{B,\mu}$. 
 As $\gamma^{-1} \in \Stab_{\Aut_X(F_n)}(N_{B,\mu})$ and $\gamma^{-1}(v_i \hat c_j v_i^{-1}) = \hat c_i \in N_{B,\mu}$, Lemma~\ref{lemma:inNBmu} implies that $r_1 \mid 1$. That is a contradiction to $r_1 > 1$. 
 
 Now let $i \in \{\mu + 1, \dots, \nu\}$ and suppose that $\gamma(\hat c_i) = v_i \hat c_j v_i^{-1}$ where $j \leq \mu$. Then $\gamma^{-1} \in \Stab_{\Aut_X(F_n)}(N_{B,\mu})$ and $ \gamma^{-1}(\hat c_j) = \gamma^{-1}(v_i^{-1}) \hat c_i \gamma^{-1}(v_i)$. But this is something that we just excluded.
 
 Elements of $C_{B,\mu}$ are either of the form $w \hat c_i^{r_\varsigma} w^{-1}$ with $w \in F_n$, $m_a(w) \in b_\varsigma$ and $i \in \{1, \dots, \mu\}$ or of the form $w \hat c_i w^{-1}$ with $w \in F_n$ and $i \in \{ \mu+1, \dots, \nu\}$. 
 For $i \in \{ \mu+1, \dots, \nu\}$ the above claim states that $\gamma(\hat c_i) \in C_{B,\mu}$ and $\gamma(w \hat c_i w^{-1}) \in C_{B,\mu}$ follows immediately. 
 If $i \in \{ 1, \dots, \mu\}$, then $h \coloneqq w \hat c_i^{r_\varsigma} w^{-1}$ with $w \in F_n$ and $m_a(w) \in b_\varsigma$ is mapped to 
 $\gamma(h) = \gamma(w) v_i \hat c_{j_i}^{r_\varsigma} v_i^{-1} \gamma(w^{-1}) \in N_{B,\mu}$. 
 Thus Lemma~\ref{lemma:inNBmu} implies that $r_\varrho \mid r_\varsigma$ where $\varrho \in \{ 1, \dots, p\}$ such that $m_a(\gamma(w) v_i) \in b_\varrho$. 
 As $\gamma(w) v_i \hat c_{j_i}^{r_\varrho} v_i^{-1} \gamma(w^{-1}) \in N_{B,\mu}$ also $\gamma^{-1}(\gamma(w) v_i \hat c_{j_i}^{r_\varrho} v_i^{-1} \gamma(w^{-1})) = w \hat c_i^{r_\varrho} w^{-1} \in N_{B,\mu}$ which implies $r_\varsigma \mid r_\varrho$. 
 Thus $\varsigma = \varrho$ and $\gamma(h) \in C_{B,\mu}$. 
 This completes 
 $\Stab_{\Aut_X(F_n)}(N_{B,\mu}) \subseteq \Stab_{\Aut_X(F_n)}(C_{B,\mu})$.
 
 The preceding paragraph also shows that for $\gamma \in \Stab_{\Aut_X(F_n)}(C_{B,\mu})$, $w \in F_n$ and $i \in \{1, \dots, \mu\}$, the elements $m_a(w)$ and $m_a(\gamma(w) v_i)$ lie in a common partition set $b_\varsigma$ of $B$. 
 Thus $m_a(v_i) = m_a(1_{F_n}) =  0$ and $m_a(w)$ and $m_a(\gamma(w))$ are in the same partition set of $B$ for all $w \in F_n$. 
 Hence it follows that $\Stab_{\Aut_X(F_n)}(C_{B,\mu}) \subseteq G_{B,\mu}$.
 
 It remains to show that $\Stab_{\Aut_X(F_n)}(C_{B,\mu}) \supseteq G_{B,\mu}$. Therefore let $\gamma \in G_{B, \mu}$ and $h \coloneqq w \hat c_i^{r_\varsigma} w^{-1} \in C_{B,\mu}$ with $w \in F_n$, $m_a(w) \in b_\varsigma$ and $i \in \{1, \dots, \mu\}$. Then $\gamma(h) = \gamma(w) v_i \hat c_{j_i}^{r_\varsigma} v_i^{-1} \gamma(w)^{-1}$ with $j_i \in \{1, \dots, \mu\}$, and $m_a( \gamma(w) v_i) = m_a(\gamma(w)) +  m_a(v_i) = m_a(\gamma(w))$ lies in the same partition set of $B$ as $m_a(w)$. Thus $\gamma(h) \in C_{B,\mu}$. 
 To that end, let $w \hat c_i w^{-1} \in C_{B,\mu}$ with $w \in F_n$ and $i \in \{ \mu+1, \dots, \nu\}$. 
 Then $\gamma(\hat c_i) = v \hat c_j v^{-1}$ with $j > \mu$ because $\gamma$ is the lift of an affine map $f$ on $X$, and this affine map permutes the singularities of $X$. 
 As the set of singularities $\{s_1, \dots, s_\mu\}$ is preserved by $f$, its complement $\{s_{\mu+1}, \dots, s_\nu\}$ is also respected by $f$. This implies $\gamma(w \hat c_i w ^{-1}) \in C_{B,\mu}$.
\end{proof}

The next lemma is the first place where we actually make use of property $(\star)$.
 
\begin{lemma}
\label{lemma:varthetaGBmu}
 Let $\vartheta \colon \Aut_X(F_n) \to \Gamma(X), \gamma_A \mapsto A$ be defined as in Section~\ref{sec:basics}. 
 Then $\vartheta(G_{B,\mu}) = \Gamma_B$.
\end{lemma}
\begin{proof}
 First let $\gamma \in G_{B,\mu}$. Then, by Definition of $G_{B,\mu}$, $m_a(w)$ and $m_a(\gamma(w))$ lie in the same partition set of $B$ for every $w \in F_n$. 
 For $A \coloneqq \vartheta(\gamma)$, arbitrary $z \in (\Z/a\Z)^{2g}$ and $w \in F_n$ such that $m_a(w) = z$ we have:
 $$\bar \varphi_a(A)(z) = \varphi_a(\gamma)(m_a(w)) = m_a(\gamma(w))   
\,.$$
 Thus $A$ respects the partition $B$ and therefore belongs to $\Gamma_B$. 
 Note that so far we did not use property $(\star)$. 
 
 Next let $A \in \Gamma_B$. 
 Recall that the affine maps $f_j$ whose derivatives $A_j$ generate $\Gamma(X)$ in property $(\star)$ (see Definition~\ref{def:star}) define via $\aff(A_j) = f_j$ a section $\aff \colon \Gamma(X) \to \Aff(\bar Y_a)$ of the map $\der$.
 
 If $A \in \Gamma_B$, then $A$ respects the partition $B$.
 Hence for every lift $\gamma$ of $A$ to $\Aut(F_n)$ and all $w \in F_n$, $m_a(w)$ and $m_a(\gamma(w))$ lie in the same partition of $B$. Now let $\gamma_A \coloneqq \lift (\aff(A))$. 
 Then by definition of the map $\aff$, for every $i \in \{1, \dots, \mu\}$, $(\aff(A))(\hat s_i) = \hat s_j$ with $j \in \{1, \dots, \mu\}$.
 Thus $\gamma_A(\hat c_i) = v_i \hat c_j v_i^{-1}$ with $m_a(v_i) = 0$. This implies that $\gamma_A \in G_{B, \mu}$. 
\end{proof}

Unfortunately the subgroup $N_{B,\mu}$ has infinite index in $\pi_1(X)$ and thereby does not define a (finite) translation surface covering $\bar X$. 
We introduce two more subsets of $N_X$ that help us to increase the subgroup $N_{B,\mu}$ in a way that guarantees that the $\vartheta$-image of the stabiliser of the enlarged subgroup remains in $\Gamma_B$.
$$ \begin{array}{rcl}
 P_X &\coloneqq& \{ v^l \mid v \in C_X, l \in \Z \}
 = \{ w {\hat c_i}^l w^{-1} \mid w \in F_n, i\in \{1,\dots,\nu\}, l \in \Z \} 
 \\
P_{B,\mu} &\coloneqq& \{ v^l \mid v \in C_{B,\mu}, l \in \Z \}\\
 &=& \{w {\hat c_i}^l w^{-1} \in P_X \mid i \in \{1, \dots, \mu\} \textrm{ and } m(w) \in b_\varsigma \Rightarrow (r_\varsigma \textrm{ divides } l) \} \\ &&\cup \; \{w {\hat c_i}^l w^{-1} \mid i \in \{\mu +1, \dots, \nu \} ,\, l \in \Z \}
\stackrel{\textrm{Lemma } \ref{lemma:inNBmu}}{=} N_{B,\mu} \cap P_X \,.
\end{array}$$
Note that because of $\gamma(C_X) = C_X$ the set $P_X$ is stabilised by every $\gamma \in \Aut_X(F_n)$.
 
\begin{cor}[see Corollary 6.9 in \cite{Sch05}]
\label{cor:PBmu}
 Let $U \leq F_n$ with $U \cap P_X = P_{B,\mu}$. 
 Then $\Stab_{\Aut_X(F_n)}(U) \subseteq G_{B,\mu}$.
\end{cor} 
\begin{proof}
 Let $\gamma \in \Stab_{\Aut_X(F_n)}(U)$.
 As $\gamma$ stabilises $U$ and $P_X$, it stabilises $P_{B,\mu} = U \cap P_X$.
 Then by $C_{B,\mu} \subseteq P_{B,\mu}$ it follows that $\gamma(C_{B,\mu}) \subseteq \gamma(P_{B,\mu}) \subseteq P_{B,\mu} \subseteq N_{B,\mu}$.
 As $N_{B,\mu} = \langle C_{B,\mu} \rangle$, this implies that $\gamma(N_{B,\mu}) \subseteq N_{B,\mu}$. 
 Thus $\gamma$ stabilises $N_{B,\mu}$. 
 By Lemma~\ref{lemma:stab=stab_mu} we have that $\gamma \in G_{B,\mu}$.
\end{proof}

Now we have everything we need in order to prove Theorem \ref{thm:GammaBAsVeechGroup}.

\begin{proof}[Proof of Theorem \ref{thm:GammaBAsVeechGroup}]
 We choose $p$ pairwise different, positive, odd numbers $r_1, \dots, r_p \in \N$, greater than $1$. 
 Then we define a translation covering $p\colon \bar Y \to \bar Y_a$ such that every preimage of a singularity $s \in \Sigma' = p_a^{-1}(p_a(S)) \subseteq \Sigma(\bar Y_a)$ with $\tilde m_a(s) \in b_\varsigma$ has ramification index $r_\varsigma$. Outside of $\Sigma'$, the covering $p$ is chosen to be unramified. 
 According to \cite{EKS84} Proposition 3.3 a topological covering with these properties always exists. 
 By lifting the translation structure from $\bar Y_a$ to $\bar Y$, we make it a translation covering.
 
 For $h = w \hat c_i w^{-1}$ with $i \in \{1, \dots, \mu\}$ and $m_a(w)\in b_\varsigma$, $h \in \pi_1(Y_a)$ is freely homotopic to a singularity $s \in \Sigma(\bar Y_a)$ with $\tilde m_a (s) \in b_\varsigma$.
 All preimages of $s$ via $p$ have ramification index $r_\varsigma$. 
 Thus $h^{r_\varsigma}$ is the smallest power of $h$ contained in $\pi_1(Y)$. 
 The ramification further implies that $w \hat c_i w^{-1} \in \pi_1(Y)$ for every $i \in \{\mu+1, \dots, \nu\}$. 
 Thus $\pi_1(Y) \cap P_X = P_{B,\mu}$ and Corollary~\ref{cor:PBmu} together with Lemma~\ref{lemma:varthetaGBmu} implies that $\Gamma(Y) \subseteq \Gamma_B$. 
 
 For the second step, we define $W \coloneqq \bigcap_{\gamma \in G_{B,\mu}} \gamma(\pi_1(Y))$. This is a finite index subgroup of $F_n$, stabilised by all $\gamma \in G_{B,\mu}$. Hence $G_{B,\mu} \subseteq \Stab_{\Aut_X(F_n)}(W)$. 
 The subgroup $W$ defines a finite translation covering $q \colon \bar Z \to \bar X$ with $\pi_1(Z) = W$, and it remains to prove $\Stab_{\Aut_X(F_n)}(\pi_1(Z)) \subseteq G_{B,\mu}$. Then, by Proposition \ref{prop:Veech group by stabiliser}, $\Gamma(Z) = \vartheta(G_{B,\mu}) = \Gamma_B$.
 
 The set $G_{B,\mu}$ is defined in a way that assures $\gamma(N_{B,\mu}) = N_{B,\mu}$ for all $\gamma \in G_{B,\mu}$. 
 Hence, for $\gamma \in G_{B,\mu}$, 
 $$P_{B,\mu} = P_X \cap N_{B,\mu} 
 = \gamma(P_X) \cap \gamma(N_{B,\mu}) 
 \stackrel{\gamma \textrm{ injective}}{=} \gamma(P_X \cap N_{B,\mu}) 
 = \gamma(P_{B,\mu})\,.$$ 
 As $P_{B,\mu} \subseteq \pi_1(Y)$ this implies $P_{B,\mu} = \bigcap_{\gamma \in G_{B,\mu}} \gamma(P_{B,\mu}) \subseteq \bigcap_{\gamma \in G_{B,\mu}} \gamma(\pi_1(Y)) = W$. 
 Of course $P_{B,\mu} \subseteq P_X$, so $P_{B,\mu} \subseteq P_X \cap W$.
 Furthermore, we have $W \subseteq \pi_1(Y)$. 
 Hence $P_X \cap W \subseteq P_X \cap \pi_1(Y) = P_{B,\mu}$. 
 Thus $P_X \cap W = P_{B,\mu}$ and consequently, by Corollary~\ref{cor:PBmu}, $\Stab_{\Aut_X(F_n)}(\pi_1(Z)) \subseteq G_{B,\mu}$.
\end{proof}

\begin{rem}
 Denote by $\mathrm{p}\Gamma(X)$ the pure Veech group of a primitive translation surface $X$. 
 Without the need of property $(\star)$ one can prove in analogy to Theorem \ref{thm:GammaBAsVeechGroup} the following result (see Theorem 2 in \cite{Fin13}):
 For every partition $B = \{b_1, \dots, b_p\}$ of $(\Z/a\Z)^{2g}$ there exists a translation covering $p \colon \bar Z \to \bar X$ with $\Gamma(Z) = \Gamma_B \cap \mathrm{p}\Gamma(X)$.
\end{rem}

\subsection{Partition stabilising groups}
\label{sec:PartStabGroups}

The groups $\Gamma_B$ are congruence groups of level $a$. 
By Theorem \ref{thm:GammaBAsVeechGroup} we are now able to find a translation covering $\bar Z \to \bar X$ for every primitive surface $\bar X$, every $a\geq 2$ and every partition $B$ of $(\Z/a\Z)^{2g}$ such that $\Gamma(Z) = \Gamma_B$  whenever $\bar X$ has property $(\star)$ in level $a$.
Obviously the next question is, which subgroups of $\Gamma(X)$ equal $\Gamma_B$ for a suitable partition $B$. 

Let $\Gamma \subseteq \Gamma(X)$ be a congruence group of level $a$. Then $\Gamma = \Gamma_B$ for the partition $B$ of $(\Z/a\Z)^{2g}$ iff $\bar{\varphi}_a (\Gamma)$ is the stabiliser of $B$ in $\bar{\varphi}_a (\Gamma(X))$. So we would like to know which congruence subgroups of $\Gamma(X)$ are stabilising groups in this sense.

A first corollary to Theorem \ref{thm:GammaBAsVeechGroup} is the following.

\begin{cor}
 \label{cor:principalCongruenceGroup}
 Let $\bar X$ be a primitive translation surface with property $(\star)$ in level $a$.
 Then the principal congruence group in $\Gamma(X)$ of level $a$ can be realised as Veech group of a translation surface. 
\end{cor}
\begin{proof}
 The principal congruence group $\Gamma(a)$ of level $a$ equals $\Gamma_B$ for 
 $$B = \{ \, \{z\} \mid z \in (\Z/a\Z)^{2g}\}\,.\qedhere$$ 
\end{proof}

Given a group $G$ acting on a set $X$ one could ask the more general question:
which subgroups of $G$ are the stabiliser of a partition $B$ of $X$? 
This question was answered in Chapter~6.5 in \cite{Sch05}:
$\bar\Gamma$ is the stabiliser of a partition iff it is the stabiliser of its orbit space (see Corollary~6.24 in \cite{Sch05}). Thus we only have to check whether $\bar{\varphi}_a (\Gamma)$ is the stabiliser of its orbit space. 

The smallest congruence group of a particular level is the principal congruence group. It is the stabiliser of its trivial orbit space. The following lemma shows that the next smallest congruence groups are stabilising groups as well.

\begin{prop}
\label{prop:index2}
 If $\bar X$ has property $(\star)$ in level $a$ and $\Gamma \leq \Gamma(X)$ is a congruence group of level $a$ with $[\Gamma: \Gamma(a)] = 2$, then $\Gamma$ is the Veech group of a covering surface of $\bar X$. 
\end{prop}
\begin{proof}
 The image of $\Gamma$ in $\Aut((\Z/a\Z)^{2g})$ has order $2$, say $\bar{\varphi}_a(\Gamma) = \{I, \bar A\}$. Thus every orbit has length $\leq 2$. 
 This implies that the image of $v \in (\Z/a\Z)^{2g}$ through $\bar A$ is uniquely determined by the orbit it lies in: if $\{v\}$ is an orbit consisting of one element, then certainly $\bar A \cdot v = v$, and if $v$ lies in $\{v, w\}$, then $\bar A \cdot v = w$. 
 
 Now suppose that there exists a $\bar C \in \bar \varphi_a(\Gamma(X)) \setminus \bar \varphi_a(\Gamma)$ that respects the orbit space of $\bar \varphi_a(\Gamma)$. 
 Then, as $\bar C \neq \bar A$ and $\bar C \neq I$, there would be $\bar \varphi_a(\Gamma)$-orbits $\{v, w\}$ and $\{v',w'\}$ such that $\bar C \cdot v = v$, $\bar C \cdot w = w$, $\bar C \cdot v' = w'$ and $\bar C \cdot w' = v'$. 
 This implies that $\bar C (v+v') = v + w'$. 
 Because of $v' \neq w'$ it follows that $v + v' \neq v + w'$ and as $\bar C$ respects the $\bar \varphi_a(\Gamma)$-orbits, $\bar A (v + v') = v+w'$. But $\bar A (v + v') = w + w' \neq v + w'$, as $w \neq v$. This is a contradiction, hence $\bar C$ does not exist.
\end{proof}

\section[Regular \texorpdfstring{$n$}{n}-gons]{Regular \texorpdfstring{\boldmath{$n$}}{n}-gons}
\label{sec:Xn}

Two classes of examples of primitive translation surfaces are the family of regular double-$n$-gons for odd $n \geq 5$ and the family of regular $2n$-gons for $n \geq 4$. 
Veech himself considered in \cite{Vee89} the family of double-$n$-gons for $n = 3$ and all $n \geq 5$. 
He constructed these surfaces by the billiard unfolding construction (see \cite{ZK75} or \cite{FK36}) in a polygon with angles $\nicefrac{\pi}{n}$, $\nicefrac{\pi}{n}$ and $\nicefrac{(n-2)\pi}{n}$. 
For even $n$, the regular double-$n$-gon is a degree two covering of the regular $n$-gon. 
We only consider the cases where the genus of the surfaces is greater than $1$, which leads to the bounds $n \geq 5$ and $n \geq 4$, respectively.
Other references concerning the double-$n$-gons are \cite{HS01} Chapter~1.7 and \cite{Vor96} Chapter~4.

\subsection[The regular double-\texorpdfstring{$n$}{n}-gon]{The regular double-\texorpdfstring{\boldmath{$n$}}{n}-gon}
\label{sec:ngon}

As the name suggests, the regular double-$n$-gon is obtained by gluing two regular $n$-gons. There is only one way to glue them that results in a translation surface: 
first identify two arbitrary sides of the two $n$-gons. 
This fixes their relative position in the plane and leads to a $2(n-1)$-gon $P$. 
Each edge in $P$ has a unique parallel edge. 
Identifying each edge with its partner gives a compact orientable surface of genus $g = \frac{n-1}{2}$ that we call $\bar X_n$. 
The translation structure of $\bar X_n$ has exactly one singular point with conical angle $(n-2) \cdot 2\pi$.
The fact that $\bar X_n$ is primitive is well known and proven e.g.\ in \cite{Fin11} Lemma~3.3. 
There, the alternative definition for a surface to be primitive (mentioned in Section \ref{sec:basics}) was used, but the given proof also meets our definition.

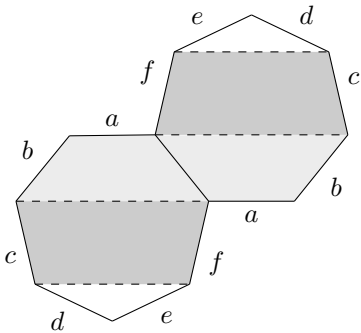
\begin{figure}[htbp]
 \centering
\begin{tikzpicture}[scale=1.3,descr/.style={fill=white,inner sep=1pt}]

\tikzset{
    myarrow/.style={->, >=latex', shorten >=1pt, thick}
} 

\coordinate (A1) at (-0.434,-0.9);
\coordinate (A2) at (0.434,-0.9);
\coordinate (A3) at (0.975,-0.223);
\coordinate (A4) at (0.782,0.623);
\coordinate (A5) at (0,1);
\coordinate (A6) at (-0.782,0.623);
\coordinate (A7) at (-0.975,-0.223);
\coordinate (A8) at  (-1.843,-0.233);
\coordinate (A9) at (-2.384,-0.9);
\coordinate (A10) at (-2.191,-1.746);
\coordinate (A11) at (-1.409,-2.123);
\coordinate (A12) at (-0.627,-1.746);

 \fill[color=gray!15] (A9) -- (A2) -- (A3) -- (A8);
 
 \fill[color=gray!40] (A7) -- (A3) -- (A4) -- (A6);
 \fill[color=gray!40] (A10) -- (A12) -- (A1) -- (A9);
 

 \draw (A1) -- node[auto,swap] {$a$} (A2); 
 \draw (0.434,-0.9) -- node[auto,swap] {$b$} (0.975,-0.223); 
 \draw (0.975,-0.223) -- node[auto,swap] {$c$} (0.782,0.623);
 \draw (0.782,0.623) -- node[auto,swap] {$d$} (0,1);
 \draw (0,1) -- node[auto,swap] {$e$} (-0.782,0.623);
 \draw (-0.782,0.623) -- node[auto,swap] {$f$} (-0.975,-0.223);
 \draw (-0.975,-0.223) -- (-0.434,-0.9);
 
 \draw (A7) -- node[auto,swap] {$a$} (A8);
 \draw (-1.843,-0.233) -- node[auto,swap] {$b$} (-2.384,-0.9);
 \draw (-2.384,-0.9) -- node[auto,swap] {$c$} (-2.191,-1.746);
 \draw (-2.191,-1.746) -- node[auto,swap] {$d$} (-1.409,-2.123);
 \draw (-1.409,-2.123) -- node[auto,swap] {$e$} (-0.627,-1.746);
 \draw (-0.627,-1.746) -- node[auto,swap] {$f$} (-0.434,-0.9);
 
 \draw[dashed] (-0.782,0.623) -- (0.782,0.623);
 \draw[dashed] (-0.975,-0.223) -- (0.975,-0.223);

 \draw[dashed] (-2.384,-0.9) -- (-0.434,-0.9);
 \draw[dashed] (-2.191,-1.746) -- (-0.627,-1.746);
 
  
\end{tikzpicture}
 \caption{Decomposition of $\bar X_7$ into horizontal cylinders.}
 \label{fig:X7cylinder}
\end{figure}

Veech determined the Veech group\index{Veech group} of $\bar X_n$ in \cite{Vee89} Theorem 5.8. It is generated by the matrices 
$$R = R(n) = \begin{pmatrix} \cos{\nicefrac{\pi}{n}} & -\sin{\nicefrac{\pi}{n}} \\ \sin{\nicefrac{\pi}{n}} & \cos{\nicefrac{\pi}{n}} \end{pmatrix} 
\; \textrm{and} \quad 
T = T(n) = \begin{pmatrix} 1 & \lambda_n \\ 0 & 1 \end{pmatrix}$$
where $\lambda_n = 2 \cot{\nicefrac{\pi}{n}}$. 
The affine map with derivative $R$ rotates the two $n$-gons and interchanges them. 
A maximal connected set of homotopic simple closed geodesics in a translation surface is called \nom{cylinder}. 
To understand the affine map $f_T$ with derivative $T$, we use the fact that $\bar X_n$ decomposes into horizontal cylinders (see Figure~\ref{fig:X7cylinder}).  
The ratio of the width to the height in each of these $\frac{n-1}{2}$ cylinders is $\lambda_n$ (see \cite{Vee89} Section~5). 
Thus the map $f_T$ is obtained by concatenating a Dehn twist of each of the cylinders.

The projective Veech group is the orientation preserving part of the Hecke triangle group with signature $(2,n,\infty)$. 
A presentation of $\Gamma(X_n)$ is $\langle R,T \mid R^{2n} = I, (T^{-1}R)^2 = R^n, R^n T = T R^n \rangle$, where $I$ is the identity element.

\begin{rem}
 As the surface $\bar X_n$ has exactly one singular point, it has property~$(\star)$ in every level. In \cite{Fin13} Section~3.4 we prove for odd $n\geq 5$ that $\Gamma(X_n) / \Gamma(2)$ is the dihedral group with $2n$ elements. For odd $n \geq 7$ every congruence group of level $2$ is the Veech group of a translation covering of $\bar X_n$. This is the result of Theorem~5 in \cite{Fin13}. Its proof uses Theorem~\ref{thm:GammaBAsVeechGroup}. 
\end{rem}

We use $\bar X_n$ as primitive base surface in Section~\ref{sec:Wohlfahrt}. There we examine the congruence levels of a congruence group in $\Gamma(X_n)$.

\subsection[The regular \texorpdfstring{$2n$}{2n}-gon]{The regular \texorpdfstring{\boldmath{$2n$}}{2n}-gon}
\label{sec:ngonEven}
The surfaces in the last section, the double-$n$-gons for odd $n \geq 5$, have exactly one singularity. The family of regular $2n$-gons with $n \geq 4$ additionally gives for odd $n$ some examples of primitive translation surfaces with more then one singularity (with two, to be precise).
These translation surfaces are obtained by identifying the parallel sides of a regular $2n$-gon. We call them $\bar X_{2n}$. 
It is well known that $\bar X_{2n}$ is primitive. 

\begin{figure}[htbp]
 \centering
 \begin{tikzpicture}[scale=1.6]

\tikzset{
    myarrow/.style={->, >=latex', shorten >=1pt, thick}
} 

 \fill (0,0) circle (1pt);

 \draw (1,0) 
 -- node[auto,swap] {$x_0$} ++(-0.191,0.588) 
 -- node[auto,swap] {$x_1$} ++(-0.5,0.363)  
 -- node[auto,swap] {$x_2$} ++(-0.618,0) 
 -- node[auto,swap] {$x_3$} ++(-0.5,-0.363) 
 -- node[auto,swap] {$x_4$} ++(-0.191,-0.588)
 -- ++(0.191,-0.588) 
 -- ++(0.5,-0.363)  
 -- ++(0.618,0) 
 -- ++(0.5,0.363) 
 -- ++(0.191,0.588);
 
 \draw[dashed] (-1,0) --(1,0);
  \draw[dashed] (-0.809,0.588) -- (0.809,0.588);
 \draw[dashed] (-0.809,-0.588) -- (0.809,-0.588);
 
 \draw[myarrow] (-0.9045,-0.294) -- (0.9045,0.294); 
 \draw[myarrow] (-0.559, -0.7695) -- (0.559, 0.7695); 
 \draw[myarrow] (0,-0.951) -- (0,0.951); 
 \draw[myarrow] (0.559,-0.7695) -- (-0.559,0.7695); 
 \draw[myarrow] (0.9045,-0.294) -- (-0.9045,0.294); 
 
  \fill (1,0) circle (1.5pt);
 \draw (0.809,0.588) node [shape=rectangle,black!40,fill,inner sep=2pt] {};
 \fill (0.309,0.951) circle (1.5pt);
 \draw (-0.309,0.951) node [shape=rectangle,black!40,fill,inner sep=2pt] {};
 \fill (-0.809,0.588) circle (1.5pt);
 \draw (-1,0) node [shape=rectangle,fill,black!40,inner sep=2pt] {};
 \fill (-0.809,-0.588) circle (1.5pt);
 \draw (-0.309,-0.951) node [shape=rectangle,black!40,fill,inner sep=2pt] {};
 \fill (0.309,-0.951) circle (1.5pt);
 \draw (0.809,-0.588) node [shape=rectangle,black!40,fill,inner sep=2pt] {};
\end{tikzpicture}
 \caption{Generators of the fundamental group of $X_{10}$.}
 \label{fig:X10fundGroup}
\end{figure}
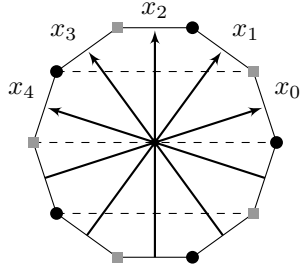

A regular $10$-gon is shown in Figure~\ref{fig:X10fundGroup}.  
The vertices of the $10$-gon are glued to two singularities, drawn as circle and rectangle in the figure. 

If $n$ is even, then $\bar X_{2n}$ has one singular point.
If $n$ is odd, then $\bar X_{2n}$ has two singular points. 
Computing the Euler characteristic we deduce that $\bar X_{2n}$ has genus
 $$g(\bar X_{2n}) = \left\{ \begin{array}{ll}
                n/2 & \textrm{, if } n \textrm{ is even}\\
		(n-1)/2 & \textrm{, if } n \textrm{ is odd}
               \end{array}
\right. \, .$$

The fundamental group of $X_{2n}$ is in both cases free in $n$ generators.
We use the following basis of the fundamental group $\pi_1(X_ {2n})$ (see Figure~\ref{fig:X10fundGroup} for an example): 
the centre of the $2n$-gon is our base point. 
Up to scaling and rotating, a regular $2n$-gon has its vertices in $(\cos (j \frac{2 \pi}{2n}), \sin(j \frac{2 \pi}{2n}))$ with $j = 0, \dots, 2n-1$. We call this normalised regular $2n$-gon $P$.
Next we number the first $n$ edges of $P$ counterclockwise with $0, \dots, n-1$, starting with the edge from $(1, 0)$ to $(\cos (\frac{2 \pi}{2n}), \sin(\frac{2 \pi}{2n}))$ and do the same with the next $n$ edges of the polygon. 
Then for every $i \in \{0, \dots, n-1\}$, there is an edge labelled with $i$ in the lower half of $P$ and a parallel edge labelled with $i$ in the upper half of $P$.
Up to homotopy, there is a unique simple closed path from the midpoint of edge $i$ in the lower half to the midpoint of edge $i$ in the upper half of $P$. 
We choose the path to cross the base point and call it $x_i$. Then the set $\{x_0, \dots, x_{n-1}\}$ is a basis of $F_n = \pi_1(X_{2n})$. 
An arbitrary element of the fundamental group $\pi_1(X_{2n})$ can be factorised in this basis by recording the labels of the crossed edges of $P$ and the directions of the crossings. This makes the basis very favourable if one defines coverings of $\bar X_{2n}$ by gluing copies of $P$.

According to Lemma~J in \cite{HS01}, the Veech group of $\bar X_{2n}$ equals the Veech group of its degree-$2$-covering, investigated by Veech. Thus Theorem~5.8 in \cite{Vee89} implies that the Veech group of the regular $2n$-gon for $n\geq 4$ is 
$\Gamma(X_{2n}) = \langle  T, R^2, R T R^{-1} \rangle$ where 
$$R = R(2n) \coloneqq \begin{pmatrix} \cos{\nicefrac{\pi}{2n}} & -\sin{\nicefrac{\pi}{2n}} \\ \sin{\nicefrac{\pi}{2n}} & \cos{\nicefrac{\pi}{2n}} \end{pmatrix} 
, \; 
T = T(2n) \coloneqq \begin{pmatrix} 1 & \lambda_{2n} \\ 0 & 1 \end{pmatrix} 
$$
and
$\lambda_{2n} = 2 \cot{\nicefrac{\pi}{2n}}$.
This is the orientation preserving part of a Hecke triangle group with signature $(n, \infty, \infty)$. 
The generators given as words in $R$ and $T$ reflect its structure as index-$2$-subgroup of the orientation preserving part of the Hecke triangle group with signature $(2, 2n , \infty)$, i.e.\ of 
$$\langle R,T \mid R^{4n} = I, (T^{-1}R)^2 = R^{2n}, R^{2n} T = T R^{2n} \rangle\,.$$
The relation $(T^{-1}R)^2 = R^{2n}$ implies that $RTR^{-1} = R^{2-2n} T^{-1}$.
Hence the third generator of $\Gamma(X_{2n})$ is redundant.
With the Reidemeister-Schreier method (see e.g.\ \cite {LS77} Chapter~II.4) and some simple transformations, one deduces the following presentation of $\Gamma(X_{2n})$:
$$\Gamma(X_{2n}) = \langle T, R^2 \mid R^{4n}, (R^2)^n T = T (R^2)^n \rangle$$
We use lifts of the affine maps with derivative $R^2$ and $T$ in $\Aut_{X_{2n}}(F_n)$ in order to determine the levels in which $\bar X_{2n}$ has property $(\star)$.
Recall that $\bar X_{2n}$ is primitive, hence there is a unique affine map $f$ with derivative $R^2$ and a unique affine map $g$ with derivative $T$ on $\bar X_{2n}$.

The affine map $f$ with derivative $R^2$ is a counterclockwise rotation of $P$ by $\nicefrac{\pi}{n}$ around its centre.  
Thus 
$$\gamma_{R^2} \colon \left\{ 
\begin{array}{lcll}
 F_n & \to & F_n \\
 x_i & \mapsto & x_{i+1} & \textrm{ for } i = 0, \dots, n-2 \\
 x_{n-1} &\mapsto& x_0^{-1}
\end{array} \right.$$
is a lift of $R^2$ to $\Aut(F_n)$.

As $|\Sigma(\bar X_{2n})| = 1$ if $n$ is even, only the surfaces $\bar X_{2n}$ for odd $n$ are examples of primitive translation surfaces with more then one singularity and add qualitatively new results to results on the double-$n$-gons. Therefore we only compute a lift of $T$ to $\Aut(F_n)$ for odd $n$. The lift for even $n$ can be obtained similarly. 

If we decompose the regular $2n$-gon into horizontal cylinders, as indicated in Figure~\ref{fig:X10fundGroup} by dashed lines, then the ratio of the width to the height of every cylinder is $\lambda_{2n}$ (see \cite{Fin11} Section~3.2). 
Hence the (unique) affine map $g$ with derivative $T$ maps each horizontal cylinder to itself and thereby shears each cylinder once. 
The image of a closed path in $X_{2n}$ under $g$ can be described as follows: every time the path traverses one of the cylinders, the appropriately oriented core curve of the cylinder is inserted into the path. Thus, in analogy to the computations in Chapter~7.3 in \cite{Fr08} for the double-$n$-gon, a lift $\gamma_T \colon F_n \to F_n$ of $T$ to $\Aut_{X_{2n}}(F_n)$ is given by 
\begin{align*}
 x_j & \mapsto \prod_{i=0}^j (x_i {x_{n-1-i}}^{-1}) \cdot x_j \cdot \prod_{i=1}^j ({x_{n-1-(j-i)}}^{-1} x_{j-i})\,,\\ 
 x_\frac{n-1}{2}  & \mapsto \prod_{i=0}^\frac{n-3}{2} (x_i {x_{n-1-i}}^{-1}) \cdot x_\frac{n-1}{2} \cdot \prod_{i=0}^\frac{n-3}{2} ({x_{n-1-(\frac{n-3}{2}-i)}}^{-1} x_{\frac{n-3}{2}-i})\,, \\
 x_{n-1-j} & \mapsto \prod_{i=0}^j (x_i {x_{n-1-i}}^{-1}) \cdot x_{n-1-j} \cdot \prod_{i=1}^j ({x_{n-1-(j-i)}}^{-1} x_{j-i})
\end{align*}
for odd $n\geq 5$ where $j = 0, \dots, \frac{n-3}{2}$.

\pagebreak

\subsection[\texorpdfstring{$\bar X_{2n}$}{X2n} and property \texorpdfstring{$(\star)$}{*}]{\texorpdfstring{\boldmath{$\bar X_{2n}$}}{X2n} and property \texorpdfstring{\boldmath{$(\star)$}}{*}}
\label{sec:X2nStar}

In this section we prove that $\bar X_{2n}$ has property $(\star)$ in level $a$ if and only if $a$ is coprime to $n$, where $n$ is odd and $n \geq 5$. 
Therefore we need to know how the affine maps with derivative $A$ act on the singularities of $\bar Y_a$ for each generator $A$ of $\Gamma(X_{2n})$. 
Proposition~\ref{prop:ActionOnSingbyHom} implies that the action on the homology helps to compute the action of the affine maps in $\bar Y_a$ on $\Sigma(\bar Y_a)$. 

The following remark recalls how the action of $\Gamma(X)$ on $H_1(\bar X, \Z/a\Z)$ is related to the action of $\Gamma(X)$ on $H_1(\bar X, \Z)$ in our context.

\begin{rem}
\label{rem:ZHom}
 Recall that the group homomorphism $\varphi_a$ which defines the action of $\Gamma(X)$ on $H_1(\bar X, \Z/a\Z)$  is the unique map with $\varphi_a(\gamma) \circ m_a = m_a \circ \gamma$ for all $\gamma \in \Aut_X(F_n)$.
The map $m_a\colon F_n \to F_n / H \cong (\Z/a\Z)^{2g}$ with $H = \nl [F_{2g},F_{2g}] \cup F_{2g}^a \cup \{c_1, \dots, c_{\nu-1}\} \nr$ factors over 
$\Z^{2g}$  through $\ab \colon F_n \to  F_n / \hat{H} \cong \Z^{2g}$, where $\hat{H} = \nl [F_{2g},F_{2g}] \cup \{c_1, \dots, c_{\nu-1}\} \nr$ and the canonical projection $\pr_a \colon \Z^{2g} \to (\Z/a\Z)^{2g}$. 
Note that the map $\pr_a \colon \Z^{2g} \to (\Z/a\Z)^{2g}$ should be interpreted as $\pr_a \colon H_1(\bar X, \Z) \to H_1(\bar X, \Z/a\Z)$.

Every $\gamma \in \Aut_X(F_n)$ descends to a unique $\ab(\gamma) \in \Aut(\Z^{2g}) \cong \GL_{2g}(\Z)$ with $\ab(\gamma) \circ \ab = \ab \circ \gamma$. 
This follows immediately by the proof of Lemma~\ref{lemma:resp_ma}. 
The map $\ab(\gamma)$ descends further to an element $\pr_a(\ab(\gamma))$ in $\Aut((\Z/a\Z)^{2g})$ and because all the descendants where unique, the following diagram commutes:
\begin{center}
\begin{tikzpicture}
   \matrix (m) [matrix of math nodes, row sep=2.5em, column sep=8em, text height=1.5ex, text depth=0.25ex]
	{ F_n & F_n\\
	  \Z^{2g} & \Z^{2g}\\
	  (\Z/a\Z)^{2g} & (\Z/a\Z)^{2g}\\};
  \path[->, font=\footnotesize]
	(m-1-1) edge node[auto] {$\gamma$} (m-1-2)
			edge [bend right=50] node[auto,swap] {$m_a$} (m-3-1)
			edge node[auto]{$\ab$} (m-2-1)
	(m-1-2) edge [bend left=50] node[auto] {$m_a$} (m-3-2)
			edge node[auto,swap]{$\ab$} (m-2-2)
	(m-2-1) edge node[auto] {$\ab(\gamma)$} (m-2-2)
			edge node[auto] {$\pr_a$} (m-3-1)
	(m-2-2) edge node[auto,swap] {$\pr_a$} (m-3-2)
	(m-3-1) edge node[auto,swap] {$\pr_a(\ab(\gamma)) = \varphi_a(\gamma)$} (m-3-2);
\end{tikzpicture}
\end{center}
\end{rem}

For the rest of this section let $n \geq 5$ be an odd natural number.
We compute $\ab(\gamma_{R^2})$ and $\ab(\gamma_T)$, where $\gamma_{R^2}$ and $\gamma_T$ are the lifts of the generators $R^2$ and $T$ of $\Gamma(X_{2n})$ in Section \ref{sec:ngonEven}.

Simple paths that are freely homotopic to the singularities $s_1$ or $s_2$, respectively, of $\bar X_{2n}$ are given by
\begin{align*}
   c_1 &= x_0 {x_1}^{-1} \cdots x_{n-3} {x_{n-2}}^{-1} x_{n-1}   \\
  \textrm{and} \quad  c_2 &= {x_0}^{-1} x_1 \cdots {x_{n-3}}^{-1} x_{n-2} {x_{n-1}}^{-1} \,.
\end{align*}
Hence the image of any subset of $\{x_0, \dots, x_{n-1}\}$ with $n-1$ elements under 
$$\ab \colon \pi_1(X_{2n}) \twoheadrightarrow \pi_1(\bar X_{2n}) / [\pi_1(\bar X_{2n}),\pi_1(\bar X_{2n})] \cong H_1(\bar X_{2n}, \Z)$$
is a basis of $H_1(\bar X_{2n}, \Z)$ and thereby induces an isomorphism $H_1(\bar X_{2n}, \Z) \cong \Z^{n-1}$. 
We choose $\{\ab(x_0), \dots, \ab(x_{\frac{n-3}{2}}), \ab(x_{\frac{n+1}{2}}), \dots, \ab(x_{n-1})\}$ as basis of $H_1(\bar X_{2n}, \Z)$. 
Let $e_j$ denote the $j$-th standard unit vector according to this basis. Then $\ab(x_j) = e_{j+1}$ for $j \in\{ 0, \dots, \frac{n-3}{2}\}$ and $\ab(x_j) = e_j$ for $j \in \{ \frac{n+1}{2}, \dots, n-1\}$. 
As $\ab(c_1) = \ab(c_2) = 0$ we obtain
$$\ab(x_\frac{n-1}{2}) = \left\{ 
\begin{array}{ll}
 \sum\limits_{i=1}^{\frac{n-1}{2}} (-1)^i \cdot e_i + \sum\limits_{i=\frac{n+1}{2}}^{n-1} (-1)^{i+1} \cdot e_i 
 & \textrm{, if } n \equiv 1 \mod 4\\ \\
 \sum\limits_{i=1}^{\frac{n-1}{2}} (-1)^{i+1} \cdot e_i + \sum\limits_{i=\frac{n+1}{2}}^{n-1} (-1)^i \cdot e_i 
 & \textrm{, if } n \equiv 3 \mod 4\\
\end{array}
\right.,$$
and hence we have: 
$$
\ab(x_\frac{n-1}{2}) =
 \begin{tikzpicture}[mymatrixenv, baseline=(current bounding box.base)]
  \matrix[mymatrix] (m)  {
	-1 \\ 1 \\ \vdots \\ -1 \\ 1 \\ 1 \\ -1 \\ \vdots \\ 1 \\ -1\\
  };
  \mymatrixbraceleft{1}{5}{$\frac{n-1}{2}$}
  \mymatrixbraceleft{6}{10}{$\frac{n-1}{2}$}
\end{tikzpicture}
 \textrm{if } n \equiv 1 \mod 4 \quad \textrm{ and } \ab(x_\frac{n-1}{2}) =
 \begin{tikzpicture}[mymatrixenv, baseline=(current bounding box.base)]
  \matrix[mymatrix] (m)  {
    1 \\ -1 \\ \vdots \\ 1 \\ -1 \\ 1 \\ 1 \\ -1 \\ 1 \\ \vdots \\ -1 \\ 1\\
  };
  \mymatrixbraceleft{1}{6}{$\frac{n-1}{2}$}
  \mymatrixbraceleft{7}{12}{$\frac{n-1}{2}$}
\end{tikzpicture}
$$
if $n \equiv 3 \mod 4$.

Define $\bar R^2 \coloneqq \ab(\gamma_{R^2})$ and $\bar T \coloneqq \ab(\gamma_{T})$. Then, with respect to the basis chosen above,
  $$\bar R^2 =
 \begin{tikzpicture}[mymatrixenv, baseline=(current bounding box.base)]
 \label{R2con1}
  \matrix[mymatrix, text height = 7pt] (m) {
    0 & \cdots &\cdots& 0 & -1 &0& \cdots&\cdots &0& -1\\
    1 & 0& \cdots &0 &1 &\vdots &&&\vdots&0\\
    0 & \ddots&\ddots&\vdots& \vdots&\vdots&&&\vdots&\vdots\\
    \vdots &\ddots&\ddots&0&-1&\vdots&&&\vdots&\vdots\\
    0 & \cdots & 0&1 & 1&0 &\cdots&\cdots&0&\vdots\\
    0&\cdots&\cdots&0&1&0&\cdots&\cdots&0&\vdots\\
    \vdots&&&\vdots&-1&1&0 &\cdots&0&\vdots\\
    \vdots&&&\vdots&\vdots&0&\ddots &\ddots&\vdots&\vdots\\
    \vdots&&&\vdots&1&\vdots&\ddots&\ddots&0&\vdots\\
    0&\cdots&\cdots&0&-1&0&\cdots&0&1&0\\
  };
  \mymatrixbraceright{2}{5}{$\frac{n-3}{2}$}
  \mymatrixbracetop{1}{4}{$\frac{n-3}{2}$}
  \mymatrixbraceleft{1}{5}{$\frac{n-1}{2}$}
  \mymatrixbraceleft{6}{10}{$\frac{n-1}{2}$}
  \mymatrixbracebottom{6}{9}{$\frac{n-3}{2}$}
\end{tikzpicture}
$$
if $n \equiv 1 \mod 4$ and
\pagebreak
$$\bar R^2 = \hspace{-0,5cm}
 \begin{tikzpicture}[mymatrixenv, baseline=(current bounding box.base)]
  \label{R2con3}
  \matrix[mymatrix,text height = 7pt] (m) {
    0 & \cdots &\cdots&\cdots& 0 & 1 &0& \cdots&\cdots&\cdots &0& -1\\
    1 & 0& \cdots &\cdots&0 &-1 &\vdots &&&&\vdots&0\\
    0 &\ddots&\ddots&&\vdots&1&\vdots&&&&\vdots&\vdots\\
    \vdots &\ddots&\ddots&\ddots&\vdots& \vdots&\vdots&&&&\vdots&\vdots\\
    \vdots &&\ddots&\ddots&0&-1&\vdots&&&&\vdots&\vdots\\
    0 & \cdots &\cdots & 0&1 & 1&0 &\cdots&\cdots&\cdots&0&\vdots\\
    0&\cdots&\cdots&\cdots&0&1&0&\cdots&\cdots&\cdots&0&\vdots\\
    \vdots&&&&\vdots&-1&1&0 &\cdots&\cdots&0&\vdots\\
    \vdots&&&&\vdots&\vdots&0&\ddots &\ddots&&\vdots&\vdots\\
    \vdots&&&&\vdots&1&\vdots & \ddots&\ddots&\ddots&\vdots&\vdots\\
    \vdots&&&&\vdots&-1&\vdots&&\ddots&\ddots&0&\vdots\\
    0&\cdots&\cdots&\cdots&0&1&0&\cdots&\cdots&0&1&0\\
  };
  \mymatrixbraceright{2}{6}{$\frac{n-3}{2}$}
  \mymatrixbracetop{1}{5}{$\frac{n-3}{2}$}
  \mymatrixbraceleft{1}{6}{$\frac{n-1}{2}$}
  \mymatrixbraceleft{7}{12}{$\frac{n-1}{2}$}
  \mymatrixbracebottom{7}{11}{$\frac{n-3}{2}$}
\end{tikzpicture}
$$ 
if $n \equiv 3 \mod 4$.

Now we compute $\bar T \cdot \ab(x_j) = \ab(\gamma_T(x_j))$ for $j \in \{0, \dots, \frac{n-3}{2}, \frac{n+1}{2}, \dots, n-1 \}$ to receive a matrix presentation for $\bar T$. For $j \in \{0, \dots, \frac{n-3}{2}\}$:
\begin{align*}
 \bar T \cdot e_{j+1} &= \ab(\gamma_T(x_j)) = \ab(\prod_{i=0}^j (x_i {x_{n-1-i}}^{-1}) \cdot x_j \cdot \prod_{i=1}^j ({x_{n-1-(j-i)}}^{-1} x_{j-i}) )\\
 &= \sum_{i=0}^j (e_{i+1} - e_{n-1-i}) + e_{j+1} + \sum_{i=1}^j (- e_{n-1-(j-i)}+ e_{j-i+1})\\
 &= \sum_{i=1}^{j+1} e_i - \sum_{i=n-1-j}^{n-1} e_i + e_{j+1} - \sum_{i=n-j}^{n-1} e_i + \sum_{i=1}^j e_i\\
 &= 2\sum_{i=1}^{j+1} e_i - e_{n-j-1} - 2 \sum_{i=n-j}^{n-1} e_i
\intertext{and}
 \bar T \cdot e_{n-1-j} &=\ab(\gamma_T(x_{n-1-j})) \\
 &= \ab(\prod_{i=0}^j (x_i {x_{n-1-i}}^{-1}) \cdot x_{n-1-j} \cdot \prod_{i=1}^j ({x_{n-1-(j-i)}}^{-1} x_{j-i}))\\
 &= \sum_{i=0}^j (e_{i+1} - e_{n-1-i}) + e_{n-1-j} + \sum_{i=1}^j (- e_{n-1-(j-i)} + e_{j-i+1})\\
 &= \sum_{i=1}^{j+1} e_i - \sum_{n-1-j}^{n-1} e_i + e_{n-1-j} - \sum_{i=n-j}^{n-1} e_i + \sum_{i=1}^{j} e_i\\
 &= 2 \sum_{i=1}^{j} e_i + e_{j+1} -2 \sum_{i=n-j}^{n-1} e_i\,.
\end{align*}
This implies
$$\bar T =
 \begin{tikzpicture}[mymatrixenv, baseline=(current bounding box.base)]
  \matrix[mymatrix] (m) {
   2 & 2 &\cdots&2&2& \cdots & 2 & 1\\
   0&\ddots&\ddots&\vdots&\vdots&\iddots&\iddots&0\\
   \vdots&\ddots&\ddots&2&2&\iddots&\iddots&\vdots\\
   0&\cdots&0&2&1&0&\cdots&\vdots\\
   0 &\cdots&0&-1&0&\cdots&\cdots&0\\
   \vdots&\iddots&\iddots&-2&-2&\ddots&&\vdots\\
   0&\iddots&\iddots&\vdots&\vdots&\ddots&\ddots&\vdots\\
   -1&-2 &\cdots &-2&-2& \cdots &-2&0\\
  };
  \mymatrixbraceleft{1}{4}{$\frac{n-1}{2}$}
  \mymatrixbraceleft{5}{8}{$\frac{n-1}{2}$}
  \mymatrixbracebottom{1}{4}{$\frac{n-1}{2}$}
  \mymatrixbracebottom{5}{8}{$\frac{n-1}{2}$}
\end{tikzpicture}.
$$
Now we return to the goal of this section and prove property $(\star)$ for $\bar X_{2n}$ in many levels.

\begin{thm}
\label{thm:X2nPropStar}
 Let $n \geq 5$ be an odd number and $a \geq 2$.
 Then the translation surface $\bar X_{2n}$ has property $(\star)$ in level $a$ if and only if $\gcd(a,n)=1$.
\end{thm}
\begin{proof}
 Let $a \geq 2$. 
 If $\gcd(a,n)=1$ then
 we determine singularities $\hat s_1$ and $\hat s_2$ in $\bar Y_a$ with $p_a(\hat s_1) = s_1$ and $p_a(\hat s_2) = s_2$ and affine maps $\hat f$ and $\hat g$ on $\bar Y_a$ with $\der(\hat f)=R^2$ and $\der(\hat g)=T$, such that $\hat f(\hat s_1) = \hat s_2$, $\hat f(\hat s_2) = \hat s_1$, $\hat g(\hat s_1) = \hat s_1$ and $\hat g(\hat s_2) = \hat s_2$. For $\gcd(a,n) > 1$ we prove that no such $\hat f$ exists.  
 
We start by investigating the affine maps $f$ and $g$ on $\bar Y_a$ with $\der(f) = R^2$ and $\der(g)= T$, given by the lifts $\gamma_{R^2}$ and $\gamma_T$ in $\Aut_X(F_n)$ from Section~\ref{sec:ngonEven}. 
Recall that 
$c_1 = x_0 {x_1}^{-1} \cdots x_{n-3} {x_{n-2}}^{-1} x_{n-1}$ and $c_2 = {x_0}^{-1} x_1 \cdots {x_{n-3}}^{-1} x_{n-2} {x_{n-1}}^{-1}$
are simple paths around the singularities of $\bar X_{2n}$.  
Obviously $\gamma_{R^2}(c_1) = x_0 c_2 {x_0}^{-1}$ and  $\gamma_{R^2}(c_2) = {x_0}^{-1} c_1 x_0$. 
A slightly longer but still easy calculation shows that $\gamma_T(c_1) = c_1$ and $\gamma_T(c_2)=c_2$.

Now we define a second lift of $R^2$ to $\Aut_X(F_n)$:
 $$\tilde \gamma_{R^2} \colon \left\{ 
\begin{array}{lcll}
 F_n & \to & F_n \\
 x_i & \mapsto & {x_0}^{-1}x_{i+1}x_0 & \textrm{ for } i = 0, \dots, n-2 \\
 x_{n-1} &\mapsto& x_0^{-1}
\end{array} \right.\,.$$
This lift fulfils $\tilde \gamma_{R^2}(c_1) = c_2$ and $\tilde \gamma_{R^2}(c_2) = {x_0}^{-2} c_1 {x_0}^2$. 
An arbitrary lift of $R^2$ to $\Aut_X(F_n)$ is given by 
$$\hat \gamma_{R^2}(w) = v \cdot \tilde \gamma_{R^2}(w) \cdot v^{-1}$$
for $v \in F_n$. 
The corresponding affine map $\hat f$ in $\Aff(\bar Y_a)$ is uniquely determined by $z \coloneqq m_a(v) \in (\Z/a\Z)^{n-1}$. 
Without loss of generality (the covering $p_a \colon \bar Y_a \to \bar X$ is normal), we choose 
$\hat c_1 \coloneqq c_1$ 
and $\hat s_1 \in \Sigma(\bar Y_a)$ such that $\hat c_1$ is freely homotopic to $\hat s_1$. 
This forces $\hat s_2 \coloneqq \hat f(\hat s_1)$ and up to conjugation with an element in $H = \pi_1(Y_a)$ it also forces
$$\hat c_2 \coloneqq \hat \gamma_{R^2}(\hat c_1) = v c_2 v^{-1}\,.$$ 
As before, we use $\{\hat c_1, \hat c_2\}$ to identify the singularities $\Sigma(\bar Y_a)$ with $\{1,2\} \times (\Z/a\Z)^{2g}$. 
By definition 
$\hat \gamma_{R^2} (\hat c_1) = \hat c_2$ and 
$$\hat \gamma_{R^2}(\hat c_2) 
= v \cdot \tilde \gamma_{R^2}( v c_2 v^{-1}) \cdot v^{-1} 
= v \cdot \tilde \gamma_{R^2}(v) \cdot {x_0}^{-2} \hat c_1 {x_0}^2 \cdot \gamma_{R^2} (v^{-1}) \cdot v^{-1} \,.$$
As $m_a(v \cdot \tilde \gamma_{R^2}(v) \cdot {x_0}^{-2}) = z + \bar R^2 z -2 e_1$, Proposition~\ref{prop:ActionOnSingbyHom} implies that 
$\hat f(1,0) = (2,0)$ and $\hat f(2,0)=(1,z + \bar R^2 z -2 e_1)$. 
Hence $\hat f$ meets the conditions of property $(\star)$ iff 
$$z + \bar R^2 z -2 e_1 = 0\,.$$

Now let $\gcd(a,n) = 1$ and hence $n$ is invertible in $\Z/a\Z$.
For $n \equiv 1 \mod 4$ we define 
$$z \coloneqq \frac{2}{n} \cdot 
\begin{pmatrix}
 \frac{n-1}{2}, & -\frac{n-3}{2}, & \dots, & 2, & -1, & 1, & -2, & \dots, & \frac{n-3}{2}, & -\frac{n-1}{2}
\end{pmatrix}^\top$$
then
$$\begin{array}{rl}
z + \bar R^2 z &= 
\frac{2}{n} 
\begin{tikzpicture}[mymatrixenv, baseline=(current bounding box.base)]
  \matrix[mymatrix,inner sep=1pt,text height = 7pt] (m)  {
    \frac{n-1}{2} \\ -\frac{n-3}{2} \\ \vdots \\ \vdots \\ 2 \\ -1 \\ 1 \\ -2 \\ \vdots \\ \vdots \\ \frac{n-3}{2} \\ -\frac{n-1}{2} \\
  };
\end{tikzpicture}
 +
\begin{tikzpicture}[mymatrixenv, baseline=(current bounding box.base)]
  \matrix[mymatrix,inner sep=1pt] (m)  {
    0 & \dots &\dots& 0 & -1 &0& \dots&\dots &0& -1\\
    1 & 0& \dots &0 &1 &\vdots &&&\vdots&0\\
    0 & \ddots&\ddots&\vdots& \vdots&\vdots&&&\vdots&\vdots\\
    \vdots &\ddots&\ddots&0&-1&\vdots&&&\vdots&\vdots\\
    0 & \dots & 0&1 & 1&0 &\dots&\dots&0&\vdots\\
    0&\dots&\dots&0&1&0&\dots&\dots&0&\vdots\\
    \vdots&&&\vdots&-1&1&0 &\dots&0&\vdots\\
    \vdots&&&\vdots&\vdots&0&\ddots &\ddots&\vdots&\vdots\\
    \vdots&&&\vdots&1&\vdots&\ddots&\ddots&0&\vdots\\
    0&\dots&\dots&0&-1&0&\dots&0&1&0\\
  };
\end{tikzpicture}
\cdot \frac{2}{n} 
\begin{tikzpicture}[mymatrixenv, baseline=(current bounding box.base)]
  \matrix[mymatrix,inner sep=1pt,text height = 7pt] (m)  {
 \frac{n-1}{2} \\ -\frac{n-3}{2} \\ \vdots \\ \vdots \\ 2 \\ -1 \\ 1 \\ -2 \\ \vdots \\ \vdots \\ \frac{n-3}{2} \\ -\frac{n-1}{2}\\
  };
\end{tikzpicture}\\ \\
 &= \frac{2}{n} 
\begin{tikzpicture}[mymatrixenv, baseline=(current bounding box.base)]
  \matrix[mymatrix,inner sep=1pt,text height = 7pt] (m)  {
  \frac{n-1}{2} + 1 + \frac{n-1}{2} \\ 
 -\frac{n-3}{2} + \frac{n-1}{2} -1 \\ 
  \vdots \\ 
 2 - 3 + 1  \\ 
 -1 + 2 - 1\\ 
 1 -1 \\ 
 -2 + 1 + 1\\  
  \vdots \\ 
 \frac{n-3}{2} - 1  - \frac{n-5}{2}\\ 
 -\frac{n-1}{2} +1 + \frac{n-3}{2}\\
  };
\end{tikzpicture}
= \frac{2}{n} 
\begin{tikzpicture}[mymatrixenv, baseline=(current bounding box.base)]
  \matrix[mymatrix,inner sep=1pt,text height = 8pt] (m)  {
  n \\ 
 0 \\ 
  \vdots \\ 
 0  \\ 
 0\\ 
 0 \\ 
 0\\ 
 \vdots \\ 
 0\\ 
 0 \\
  };
\end{tikzpicture}
= 2 \cdot e_1\,.
\end{array}
$$
Analogously one checks for $n \equiv 3 \mod 4$ and 
$$z \coloneqq \frac{2}{n} \cdot
\begin{pmatrix}
 \frac{n-1}{2}, & -\frac{n-3}{2}, & \dots, & -2, & 1, & -1, & 2, & \dots, & \frac{n-3}{2}, & -\frac{n-1}{2}
\end{pmatrix}^\top$$ 
that $z + \bar R^2 z = 2 \cdot e_1$.

As $\hat c_1 = c_1$, $g(\hat s_1) = \hat s_1$. Thus there are affine maps on $\bar Y_a$, fulfilling the conditions of property $(\star)$ iff $g$ also satisfies  $g(\hat s_2) = \hat s_2$, i.e.\ iff $g(2,0) = (2,0)$.
As
$$\gamma_T(\hat c_2) = \gamma_T(v c_2 v^{-1}) = \gamma_T(v) c_2 \gamma_T(v)^{-1} =  \gamma_T(v) v^{-1} \hat c_2 v \gamma_T(v)^{-1},$$
$g(2,0) = (2,0)$ iff  $m_a(\gamma_T(v) v^{-1}) = \bar T z -z = 0$.

For $n \equiv 1 \mod 4$ we have 
$$\begin{array}{l}
\bar T z 
=
\begin{tikzpicture}[mymatrixenv, baseline=(current bounding box.base)]
  \matrix[mymatrix] (m)  {
 2 & 2 &\dots&\dots&\dots& \dots & 2 & 1\\
 0&\ddots&\ddots&&&\iddots&\iddots&0\\
   \vdots&\ddots&\ddots&2&2&\iddots&\iddots&\vdots\\
   \vdots&&0&2&1&\iddots&&\vdots\\
   \vdots &&0&-1&0&&&\vdots\\
   \vdots&\iddots&\iddots&-2&-2&\ddots&&\vdots\\
   0&\iddots&\iddots&&&\ddots&\ddots&\vdots\\
   -1&-2 &\dots &\dots&\dots& \dots &-2&0\\
  };
\end{tikzpicture}
\cdot \frac{2}{n}
\begin{tikzpicture}[mymatrixenv, baseline=(current bounding box.base)]
  \matrix[mymatrix] (m)  {
 \frac{n-1}{2} \\ -\frac{n-3}{2} \\ \vdots \\ 2 \\ -1 \\ 1 \\ -2 \\ \vdots \\ \frac{n-3}{2} \\ -\frac{n-1}{2}\\
  };
\end{tikzpicture}\\
\quad = \frac{2}{n}
\begin{tikzpicture}[mymatrixenv, baseline=(current bounding box.base)]
  \matrix[mymatrix] (m)  {
  n-1 + 2 \sum_{i=1}^{\frac{n-3}{2}} (i-i) - \frac{n-1}{2} \\ 
  -(n-3) + 2 \sum_{i=1}^{\frac{n-5}{2}} (i-i) + \frac{n-3}{2} \\ 
  \vdots \\ 
  4-2+2-2 \\ -2+1 \\ 1 \\ -2+2-2 \\ \vdots \\ \frac{n-3}{2} + 2 \sum_{i=1}^{\frac{n-5}{2}}(i-i) \\ -\frac{n-1}{2} + 2 \sum_{i=1}^{\frac{n-3}{2}}(i-i)\\
  };
\end{tikzpicture}
= z
\end{array}$$
and similarly for $n \equiv 3 \mod 4$ one checks that $\bar T z = z$.

 Now let $\gcd(a,n) > 1$. We saw at the beginning of this proof that there is an affine map $\hat f$ on $\bar Y_a$ with $\der(\hat f) = R^2$ and $\hat s_1, \hat s_2 \in \Sigma(\bar Y_a)$ such that $\hat f(\hat s_1) = \hat s_2$ and $\hat f(\hat s_2) = \hat s_1$ iff there exists a $z \in (\Z/a\Z)^{n-1}$ such that $z + \bar R^2 z = 2 e_1$. 
 
 For $n \equiv 1 \mod 4$, this gives the following system of $n-1$ linear equations in $z = (z_1, \dots, z_{n-1})^\top$ over $\Z/a\Z$ (see $\bar R^2$ on page~\pageref{R2con1}):
 $$\begin{array}{lrcll}
    \textrm{Equation } 1: &  2 &=& z_1 - z_\frac{n-1}{2} - z_{n-1} \\
    \textrm{Equation } i: &  0 &=& z_{i-1} + z_i + (-1)^i \, z_\frac{n-1}{2} & \textrm{for } i = 2, \dots, \frac{n-3}{2} \\
    \textrm{Equation } \frac{n-1}{2}: & 0 &=& z_\frac{n-3}{2} +2 \, z_\frac{n-1}{2}\\
    \textrm{Equation } \frac{n+1}{2}: & 0 &=& z_\frac{n-1}{2} + z_\frac{n+1}{2}\\
    \textrm{Equation } i: & 0 &=& (-1)^{i+1} \, z_\frac{n-1}{2} + z_{i-1} + z_i  & \textrm{for } i = \frac{n+3}{2}, \dots, n-1 
   \end{array}
 $$
 If we add up $(-1)^i$-times the $i$-th equation for $i = 1, \dots, \frac{n-1}{2}$ and $(-1)^{i+1}$-times the $i$-th equation for $i = \frac{n+1}{2}, \dots, n-1$ we get:
 \begin{align*}
     -2 &= - (z_1 - z_\frac{n-1}{2} - z_{n-1}) 
 + \sum_{i=2}^\frac{n-3}{2} (-1)^i (z_{i-1} + z_i + (-1)^i \, z_\frac{n-1}{2})
 + z_\frac{n-3}{2}\\
 &\quad +2 \, z_\frac{n-1}{2} + z_\frac{n-1}{2} + z_\frac{n+1}{2}
 + \sum_{\frac{n+3}{2}}^{n-1} (-1)^{i+1} ( (-1)^{i+1} z_\frac{n-1}{2} + z_{i-1} + z_i )\\
 &= n \, z_\frac{n-1}{2}
 \end{align*}
 As $n$ is odd $\gcd(a,n) > 1$ implies that $b \coloneqq \gcd(a,n) > 2$. So $-2 = n \, z_\frac{n-1}{2}$ implies $-2 \frac{a}{b} = n \frac{a}{b} \, z_\frac{n-1}{2} = a \frac{n}{b} \, z_\frac{n-1}{2} \equiv 0 \mod a$. But $-2 \frac{a}{b} \not\equiv 0 \mod a$ as $b >2$. Hence there is no solution to $z + \bar R^2 z = 2 e_1$ if $\gcd(a,n) > 1$.
 
 For $n \equiv 3 \mod 4$, we get a very similar system of $n-1$ linear equations in $z = (z_1, \dots, z_{n-1})^\top$ over $\Z/a\Z$ (see $\bar R^2$ on page~\pageref{R2con3}):
  $$\begin{array}{lrcll}
    \textrm{Equation } 1: &  2 &=& z_1 + z_\frac{n-1}{2} - z_{n-1} \\
    \textrm{Equation } i: &  0 &=& z_{i-1} + z_i + (-1)^{i+1} \, z_\frac{n-1}{2} & \textrm{for } i = 2, \dots, \frac{n-3}{2} \\
    \textrm{Equation } \frac{n-1}{2}: & 0 &=& z_\frac{n-3}{2} +2 \, z_\frac{n-1}{2}\\
    \textrm{Equation } \frac{n+1}{2}: & 0 &=& z_\frac{n-1}{2} + z_\frac{n+1}{2}\\
    \textrm{Equation } i: & 0 &=& (-1)^i \, z_\frac{n-1}{2} + z_{i-1} + z_i  & \textrm{for } i = \frac{n+3}{2}, \dots, n-1 
   \end{array}
 $$
 Here we add up $(-1)^{i+1}$-times the $i$-th equation for $i = 1, \dots, \frac{n-1}{2}$ and $(-1)^{i}$-times the $i$-th equation for $i = \frac{n+1}{2}, \dots, n-1$:
  \begin{align*}
     2 &= z_1 + z_\frac{n-1}{2} - z_{n-1}
 + \sum_{i=2}^\frac{n-3}{2} (-1)^{i+1} (z_{i-1} + z_i + (-1)^{i+1} \, z_\frac{n-1}{2})
 + z_\frac{n-3}{2} +2 \, z_\frac{n-1}{2}\\
 &\quad + z_\frac{n-1}{2} + z_\frac{n+1}{2}
 + \sum_{\frac{n+3}{2}}^{n-1} (-1)^i ( (-1)^i \, z_\frac{n-1}{2} + z_{i-1} + z_i)\\
 &= n \, z_\frac{n-1}{2}
 \end{align*}
 As above, this implies $2 \frac{a}{b} = a \frac{n}{b} \, z_\frac{n-1}{2} \equiv 0 \mod a$. But $2 \frac{a}{b} \not\equiv 0 \mod a$ as $b >2$.

Thus there is no $f \in \Aff(\bar Y_a)$ with $\der(f) = R^2$ such that $f$ fixes a subset of $\Sigma(\bar Y_a)$ of cardinality $2$. 
It follows that $\bar X_{2n}$ does not have property $(\star)$ in level $a$ if $\gcd(a,n) > 1$.
\end{proof}

\section{Congruence levels}
\label{sec:Wohlfahrt}

Theorem \ref{thm:GammaBAsVeechGroup} is a strong motivation to investigate congruence groups. 
We want to determine the congruence levels of a congruence group and we would like to know how its various congruence levels are related. 
It is however very difficult for a given finite index subgroup to decide whether it is a congruence group and if so for which congruence levels. 
For finite index subgroups $\Gamma$ of the Veech group of the regular double-$n$-gons $\bar X_n$ for odd $n \geq 5$, we show how one can restrict the candidates for the potential congruence levels of $\Gamma$. We additionally give a sufficient condition for $\Gamma \leq \Gamma(X_n)$ to be not a congruence group and an example of a non-congruence group in $\Gamma(X_5)$.

A first general observation on different levels of a congruence group follows easily by Remark~\ref{rem:ZHom}:

Let $\bar X$ be a primitive translation surface, and $a,k \geq 1$. 
Then $$\Gamma(ka) \subseteq \Gamma(a) \subseteq \Gamma(X)\,.$$

This immediately implies that every congruence group of level $a$ in $\Gamma(X)$ is also a congruence group of level $k \cdot a$  for every $k \geq 1$.
Thus every congruence group has multiple congruence levels. 

\begin{defn}
 If $\Gamma \leq \Gamma(X)$ is a congruence group of level $a$, then $a$ is called \nom{minimal congruence level} of $\Gamma$, if $\Gamma$ is not a congruence group of a level $b$ with $b \mid a$.
\end{defn}

Note that if $\bar X$ is the once-punctured torus then $\Gamma(X) = \SL_2(\Z)$ and every congruence group has a unique minimal congruence level, whereas this is not clear for congruence groups in the Veech group of other primitive translation surfaces.

A second consequence of Remark~\ref{rem:ZHom} and the Chinese Remainder Theorem is that for $a,b \in \N$ with $\gcd(a,b)=1$ we have $$\Gamma(a) \cap \Gamma(b) = \Gamma(ab)\,.$$ 

\subsection{Parabolic elements in \texorpdfstring{\boldmath{$\Gamma(X_n)$}}{Xn}}
Let $n \geq 5$ be odd. We investigate the regular double-$n$-gon, $\bar X_n$, as introduced in Section \ref{sec:ngon}. 
In the following section we generalise the Wohlfahrt level to subgroups of $\Gamma(X_n)$. 
This level definition depends on the parabolic elements in the Veech group. 
Thus the goal of this section is to find the parabolic elements in the principal congruence groups $\Gamma(a) \subseteq \Gamma(X_n)$.

\begin{figure}[htbp]
 \centering
 \begin{center}
\begin{tikzpicture}[scale=1.3,descr/.style={fill=white,inner sep=1pt}]

\tikzset{
    myarrow/.style={->, >=latex', shorten >=1pt, thick}
} 

\coordinate (A1) at (-0.434,-0.9);
\coordinate (A2) at (0.434,-0.9);
\coordinate (A3) at (0.975,-0.223);
\coordinate (A4) at (0.782,0.623);
\coordinate (A5) at (0,1);
\coordinate (A6) at (-0.782,0.623);
\coordinate (A7) at (-0.975,-0.223);
\coordinate (A8) at  (-1.843,-0.233);
\coordinate (A9) at (-2.384,-0.9);
\coordinate (A10) at (-2.191,-1.746);
\coordinate (A11) at (-1.409,-2.123);
\coordinate (A12) at (-0.627,-1.746);

 \fill[color=gray!15] (A9) -- (A2) -- (A3) -- (A8);
 
 \fill[color=gray!40] (A7) -- (A3) -- (A4) -- (A6);
 \fill[color=gray!40] (A10) -- (A12) -- (A1) -- (A9);
 

 
 \draw[myarrow] (-2.1135,-0.5615) -- node[pos=0.3,fill=gray!15] {$x_1$} (0.7045,-0.5615); 
 \draw[myarrow] (-0.8785,0.2) -- node[midway,fill=gray!40] {$x_2$} (0.8785,0.2); 
 \draw[myarrow] (-2.2875,-1.323) -- node[midway,fill=gray!40] {$x_2$} (-0.5305,-1.323); 
 \draw[myarrow] (-0.391,0.8115) -- node[midway,fill=white] {$x_3$} (0.391,0.8115); 
 \draw[myarrow] (-1.8,-1.9345) -- node[midway,fill=white] {$x_3$} (-1.018,-1.9345); 

 \draw (-0.434,-0.9) -- node[auto,swap] {$a$} (0.434,-0.9); 
 \draw (0.434,-0.9) -- node[auto,swap] {$b$} (0.975,-0.223); 
 \draw (0.975,-0.223) -- node[auto,swap] {$c$} (0.782,0.623);
 \draw (0.782,0.623) -- node[auto,swap] {$d$} (0,1);
 \draw (0,1) -- node[auto,swap] {$e$} (-0.782,0.623);
 \draw (-0.782,0.623) -- node[auto,swap] {$f$} (-0.975,-0.223);
 \draw (-0.975,-0.223) -- (-0.434,-0.9);
 
 \draw (-0.975,-0.223) -- node[auto,swap] {$a$} (-1.843,-0.233);
 \draw (-1.843,-0.233) -- node[auto,swap] {$b$} (-2.384,-0.9);
 \draw (-2.384,-0.9) -- node[auto,swap] {$c$} (-2.191,-1.746);
 \draw (-2.191,-1.746) -- node[auto,swap] {$d$} (-1.409,-2.123);
 \draw (-1.409,-2.123) -- node[auto,swap] {$e$} (-0.627,-1.746);
 \draw (-0.627,-1.746) -- node[auto,swap] {$f$} (-0.434,-0.9);

 \draw[dashed] (-0.975,-0.223) -- (0.975,-0.223);
 \draw[dashed] (-0.782,0.623) -- (0.782,0.623);
 \draw[dashed] (-2.384,-0.9) -- (-0.434,-0.9);
 \draw[dashed] (-2.191,-1.746) -- (-0.627,-1.746);
 
   \draw (2,-0.5615) node {cylinder $1$};
  \draw (2,0.2) node {cylinder $2$};
  \draw (2,0.8115) node {cylinder $3$};
  \draw (-3.3,-0.5615) node {cylinder $1$};
  \draw (-3.3,-1.323) node {cylinder $2$};
  \draw (-3.3,-1.9345) node {cylinder $3$};

\end{tikzpicture}
\end{center}
\caption{Core curves of the horizontal cylinders in $\bar X_7$.}
\label{fig:coreCurves}
\end{figure}
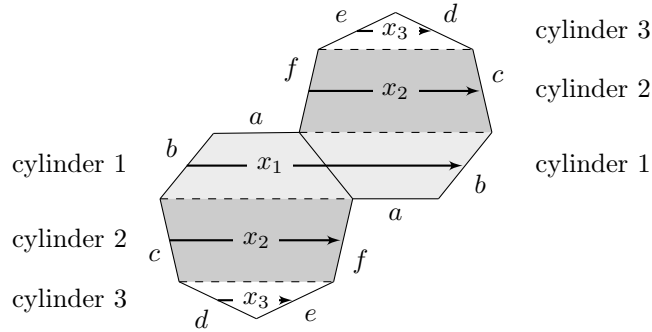

Recall that the regular double-$n$-gon $\bar X_n$ for odd $n \geq 5$ decomposes into horizontal cylinders and that the parabolic generator $T$ of the Veech group of $\bar X_n$ is the parabolic element associated to this cylinder decomposition. So $T$ is a single Dehn twist on every horizontal cylinder.

As shown in Figure~\ref{fig:coreCurves}, we arrange the two regular $n$-gons of $X_n$ such that they both have a horizontal edge. 
The right one lies above its horizontal edge and the left one lies below its horizontal edge. 
We number the cylinders in the right polygon from the bottom up with $1$ to $\frac{n-1}{2}$. 

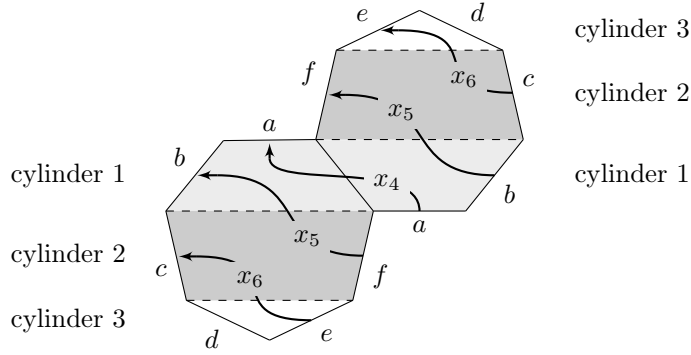
\begin{figure}[hbtp]
 \centering
%
%

\begin{center}
\begin{tikzpicture}[scale=1.4,descr/.style={fill=white,inner sep=1pt}]

\tikzset{
    myarrow/.style={->, >=latex', shorten >=1pt, thick}
} 

\coordinate (A1) at (-0.434,-0.9);
\coordinate (A2) at (0.434,-0.9);
\coordinate (A3) at (0.975,-0.223);
\coordinate (A4) at (0.782,0.623);
\coordinate (A5) at (0,1);
\coordinate (A6) at (-0.782,0.623);
\coordinate (A7) at (-0.975,-0.223);
\coordinate (A8) at  (-1.843,-0.233);
\coordinate (A9) at (-2.384,-0.9);
\coordinate (A10) at (-2.191,-1.746);
\coordinate (A11) at (-1.409,-2.123);
\coordinate (A12) at (-0.627,-1.746);

 \fill[color=gray!15] (A9) -- (A2) -- (A3) -- (A8);
 
 \fill[color=gray!40] (A7) -- (A3) -- (A4) -- (A6);
 \fill[color=gray!40] (A10) -- (A12) -- (A1) -- (A9);
 

 \draw[thick] (-1.409,-0.4) .. controls +(-90:0.3cm) and +(90:0.5cm) .. (0,-0.9) node[pos=0.7,fill=gray!15] {$x_4$};
 \draw[myarrow] (-1.409,-0.4)  -- (-1.409,-0.233);
 
 \draw[myarrow] (-0.7,0.2)  -- (-0.8785,0.2);
 \draw[thick] (-0.7,0.2) .. controls +(0:0.9cm) and +(180:0.9cm) .. (0.7045,-0.5615) node[pos=0.3,fill=gray!40] {$x_5$};
 \draw[myarrow] (-1.9,-0.5615)  -- (-2.1135,-0.5615);
 \draw[thick] (-1.9,-0.5615) .. controls +(0:0.9cm) and +(180:0.9cm) .. (-0.5305,-1.323) node[pos=0.7,fill=gray!40] {$x_5$};
 
 \draw[myarrow] (-0.2,0.8115) -- (-0.391,0.8115);
 \draw[thick] (-0.2,0.8115) .. controls +(0:0.9cm) and +(-180:0.9cm) .. (0.8785,0.22)
 node[pos=0.7,fill=gray!40] {$x_6$};
 \draw[myarrow] (-2.1,-1.33) -- (-2.2875,-1.33);
 \draw[thick] (-2.1,-1.33) .. controls +(0:0.9cm) and +(-180:0.9cm) .. (-1.018,-1.9345)
 node[pos=0.4,fill=gray!40] {$x_6$};

 \draw (-0.434,-0.9) -- node[auto,swap] {$a$} (0.434,-0.9); 
 \draw (0.434,-0.9) -- node[auto,swap] {$b$} (0.975,-0.223); 
 \draw (0.975,-0.223) -- node[auto,swap] {$c$} (0.782,0.623);
 \draw (0.782,0.623) -- node[auto,swap] {$d$} (0,1);
 \draw (0,1) -- node[auto,swap] {$e$} (-0.782,0.623);
 \draw (-0.782,0.623) -- node[auto,swap] {$f$} (-0.975,-0.223);
 \draw (-0.975,-0.223) -- (-0.434,-0.9);
 
 \draw (-0.975,-0.223) -- node[auto,swap] {$a$} (-1.843,-0.233);
 \draw (-1.843,-0.233) -- node[auto,swap] {$b$} (-2.384,-0.9);
 \draw (-2.384,-0.9) -- node[auto,swap] {$c$} (-2.191,-1.746);
 \draw (-2.191,-1.746) -- node[auto,swap] {$d$} (-1.409,-2.123);
 \draw (-1.409,-2.123) -- node[auto,swap] {$e$} (-0.627,-1.746);
 \draw (-0.627,-1.746) -- node[auto,swap] {$f$} (-0.434,-0.9);

 \draw[dashed] (-0.975,-0.223) -- (0.975,-0.223);
 \draw[dashed] (-0.782,0.623) -- (0.782,0.623);
 \draw[dashed] (-2.384,-0.9) -- (-0.434,-0.9);
 \draw[dashed] (-2.191,-1.746) -- (-0.627,-1.746);
 
  \draw (2,-0.5615) node {cylinder $1$};
  \draw (2,0.2) node {cylinder $2$};
  \draw (2,0.8115) node {cylinder $3$};
  \draw (-3.3,-0.5615) node {cylinder $1$};
  \draw (-3.3,-1.323) node {cylinder $2$};
  \draw (-3.3,-1.9345) node {cylinder $3$};

\end{tikzpicture}
\end{center}

\caption{Curves connecting neighbouring cylinders in $\bar X_7$.}
\label{fig:X7fundGroup}
\end{figure}

The fundamental group of $X_n$ is free of rank $n-1$. 
We amend the core curves $x_1, \dots, x_{\frac{n-1}{2}}$ of the horizontal cylinders to a basis of the fundamental group. 
Generator $x_\frac{n+1}{2}$ is a simple path from the bottom of cylinder $1$ in the right polygon to the top of cylinder $1$ in the left polygon.
For $i \in \{2, \dots, \frac{n-1}{2}\}$ the generators $x_{\frac{n-1}{2} +i}$ are chosen to connect neighbouring cylinders, as is shown in Figure~\ref{fig:X7fundGroup} for $\bar X_7$.

\begin{prop}
\label{prop:ParabInXn}
Let $a \geq 2$ and let $\Gamma(a)$ be the principal congruence group of level $a$ in $\Gamma(X_n)$. For $b \geq 1$: 
 $$T^b \in \Gamma(a) \Leftrightarrow a \mid b$$
\end{prop}
\begin{proof}
 At first we want to understand, what $T^b$ does to an arbitrary element $v$ in the fundamental group:
 each time $v$ crosses the cylinder $i$ from bottom to top, $T^b v$ follows the core curve $b$ times around the waist of the cylinder, i.e.\ ${x_i}^b$ is inserted into $v$. 
 If $v$ crosses cylinder $i$ from top to bottom, ${x_i}^{-b}$ is inserted at this particular position (see Figure~\ref{fig:Tv}).  

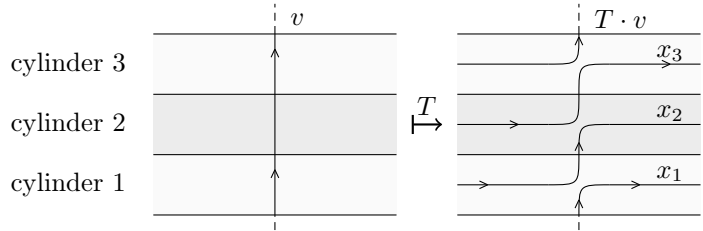
\begin{figure}[htbp]
 \centering
 \begin{tikzpicture}[scale=0.4]

  \fill[color=gray!4] (-4,4) rectangle (4,6);
  \fill[color=gray!15] (-4,2) rectangle (4,4);
  \fill[color=gray!4] (-4,0) rectangle (4,2);
 
  \draw[thin] (-4,6) -- (4,6); 
  \draw[thin] (-4,4) -- (4,4);
  \draw[thin] (-4,2) -- (4,2);
  \draw[thin] (-4,0) -- (4,0);
  
  \draw[densely dashed] (0,-0.5) --(0,0);
  \draw (0,0) -- (0,6);
  \draw[densely dashed] (0,6) --(0,7);
  \draw (0,1.5) -- +(295:0.3);
  \draw (0,1.5) -- +(245:0.3);
  \draw (0,5.5) -- +(295:0.3);
  \draw (0,5.5) -- +(245:0.3);
  
  \draw (0.2,6.5) node[right] {$v$};
  
  \draw (-9,1) node[right] {cylinder $1$};
  \draw (-9,3) node[right] {cylinder $2$};
  \draw (-9,5) node[right] {cylinder $3$};
  
 
  \begin{scope}[xshift=10cm] 
  \draw[|->,thick] (-5.5,3) -- (-4.5,3);
  \draw (-5,3) node[above] {$T$};
   
  \fill[color=gray!4] (-4,4) rectangle (4,6);
  \fill[color=gray!15] (-4,2) rectangle (4,4);
  \fill[color=gray!4] (-4,0) rectangle (4,2);
 
  \draw[thin] (-4,6) -- (4,6); 
  \draw[thin] (-4,4) -- (4,4);
  \draw[thin] (-4,2) -- (4,2);
  \draw[thin] (-4,0) -- (4,0);

  \draw (0.2,6.5) node[right] {$T\cdot v$};
  
  \draw[densely dashed] (0,-0.5) -- (0,0);
  
  \draw (0,0)
	.. controls (0,1) and (0,1) .. (1,1);
  \draw (0,0.5) -- +(295:0.3);
  \draw (0,0.5) -- +(245:0.3);
  \draw (1,1) -- (4,1);
  \draw (3,0.8) node[above] {$x_1$};
  \draw (2,1) -- +(155:0.3);
  \draw (2,1) -- +(205:0.3);
  
  \draw (-4,1) -- (-1,1);
  \draw (-3,1) -- +(155:0.3);
  \draw (-3,1) -- +(205:0.3);
  \draw (-1,1)
	.. controls (0,1) and (0,1) .. (0,2);
	
  \draw (0,2)
	.. controls (0,3) and (0,3) .. (1,3);
  \draw (0,2.4) -- +(295:0.3);
  \draw (0,2.4) -- +(245:0.3);
  
  \draw (1,3) -- (4,3);
  \draw (3,2.8) node[above] {$x_2$};
  \draw (-4,3) -- (-1,3);
  \draw (-2,3) -- +(155:0.3);
  \draw (-2,3) -- +(205:0.3);
  \draw (-1,3)
	.. controls (0,3) and (0,3) .. (0,4);
	
  \draw (0,4)
	.. controls (0,5) and (0,5) .. (1,5);
  \draw (1,5) -- (4,5);
  \draw (3,5) -- +(155:0.3);
  \draw (3,5) -- +(205:0.3);
  \draw (3,4.8) node[above] {$x_3$};
  
  \draw (-4,5) -- (-1,5);
  \draw (-1,5)
	.. controls (0,5) and (0,5) .. (0,6);
  \draw (0,5.9) -- +(295:0.3);
  \draw (0,5.9) -- +(245:0.3);
  
  \draw[densely dashed] (0,6) -- (0,7);
  \end{scope}
\end{tikzpicture}
\caption{Application of $T$ to $v$.}
\label{fig:Tv}
\end{figure}

Next we decompose $T^b v$ into elements of the fundamental group (see Figure \ref{fig:decompositionTv}).  
Therefore we follow $v$ until we reach the first insertion of a ${x_i}^{\pm b}$, then follow ${x_i}^{\pm b}$ and afterwards walk back to the base point along $v^{-1}$. This is a closed path, i.e.\ an element in the fundamental group of $X_n$. 
Next we follow $v$ until the second insertion of a ${x_i}^{\pm b}$, walk along ${x_i}^{\pm b}$ and back to the base point along $v^{-1}$. 
We continue until there is no further insertion of a ${x_i}^{\pm b}$ and finally follow $v$ entirely. 
Up to homotopy, i.e.\ the walking along $v$ back and forth, we walked exactly once along $T^b v$. 
This proves that we can decompose $T^b v$ as $T^b v = \prod_{j=1}^k (v_j x_{i_j} {v_j}^{-1})^{\pm b} \cdot v$ into elements of $\pi_1(X)$, where $i_j \in \{1, \dots, \frac{n-1}{2}\}$, $k \in \N$ and the $v_j$ are sub-paths of $v$. 

\begin{figure}[htbp]
 \centering 
 \begin{tikzpicture}[scale=0.45]
  
   
  \fill[color=gray!4] (-4,4) rectangle (4,6);
  \fill[color=gray!15] (-4,2) rectangle (4,4);
  \fill[color=gray!4] (-4,0) rectangle (4,2);
 
  \draw[thin] (-4,6) -- (4,6); 
  \draw[thin] (-4,4) -- (4,4);
  \draw[thin] (-4,2) -- (4,2);
  \draw[thin] (-4,0) -- (4,0);
  
  \draw[densely dashed] (0,-1) -- (0,0);
  \draw (0,-0.3) -- +(295:0.2);
  \draw (0,-0.3) -- +(245:0.2);

  \draw (0,0) .. controls +(50:1cm) and +(180:1cm) .. (1.5,1);

  \draw (1.5,1) -- (4,1);
  \draw (3,0.8) node[above] {$x_1$};
  \draw (2,1) -- +(155:0.2);
  \draw (2,1) -- +(205:0.2);
  
  \draw (-4,1) -- (-1.5,1);
  \draw (-3,1) -- +(155:0.2);
  \draw (-3,1) -- +(205:0.2);
  \draw (-1.5,1) .. controls +(0:1cm) and +(130:1cm) .. (0,0);

  \draw (0,0) .. controls +(65:1cm) and +(180:1cm) .. (2,3);
  \draw (0.62,1.5) -- +(275:0.2);
  \draw (0.62,1.5) -- +(225:0.2);
  \draw (2.5,3) -- +(155:0.2);
  \draw (2.5,3) -- +(205:0.2);
  
  \draw (2,3) -- (4,3);
  \draw (3,2.8) node[above] {$x_2$};
  \draw (-4,3) -- (-2,3);
  \draw (-2,3) -- +(155:0.2);
  \draw (-2,3) -- +(205:0.2);
  \draw (-2,3) .. controls +(0:1cm) and +(115:1cm) .. (0,0);
  \draw (-0.715,1.7) -- +(142:0.2);
  \draw (-0.715,1.7) -- +(92:0.2);

  \draw (0,0) .. controls +(80:1cm) and +(-180:1cm) .. (2,5);
  \draw (2,5) -- (4,5);
  \draw (3,5) -- +(155:0.2);
  \draw (3,5) -- +(205:0.2);
 
  \draw (-0.98,3.7) -- +(137:0.2);
  \draw (-0.98,3.7) -- +(87:0.2); 

  \draw (0,3) -- +(295:0.2);
  \draw (0,3) -- +(245:0.2);
  \draw (3,4.9) node[above] {$x_3$};
  
   \draw (-4,5) -- (-2,5);
   \draw (-2,5) .. controls +(0:1cm) and +(100:1cm) .. (0,0);

  \draw (0,0) -- (0,6);
  
  \draw[densely dashed] (0,6) -- (0,7);
  \draw (0,6.3) -- +(295:0.2);
  \draw (0,6.3) -- +(245:0.2);
\end{tikzpicture}
\caption{Decomposition of $T \cdot v$.}
\label{fig:decompositionTv}
\end{figure}
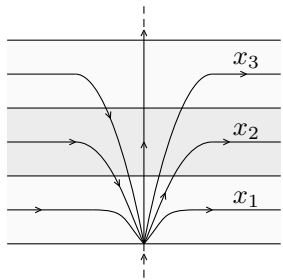

Let $\tilde{w}$ denote the image of $w \in \pi_1(X_n)$ in the absolute homology $H_1(\bar X_n, \Z/a\Z)$. Then 
$$T^b \tilde v = b \cdot \tilde{h} + \tilde v \in H_1(\bar X_n,\Z/a\Z)$$
where $h = \prod_{j=1}^k v_j  {x_{i_j}}^{\pm 1} {v_j}^{-1}$.
Now it is obvious that $a \mid b$ is a sufficient condition for $T^b$ to act trivially on the homology. 

It remains to show that $T^b \in \Gamma(a)$ implies that $a$ divides $b$.
The image of $x_\frac{n+1}{2}$ under $T^b$ is $x_\frac{n+1}{2} {x_1}^b$.
Thus $T^b \cdot \tilde x_\frac{n+1}{2} = \tilde x_\frac{n+1}{2} + b \cdot \tilde x_1$.
As $\bar X_n$ has exactly one singular point, the generators of $\pi_1(X_n)$ form a basis $\{ \tilde x_1, \dots, \tilde x_{n-1}\}$ of the homology of the closed surface $H_1(\bar X_n, \Z/a\Z)$. 
Thus, if $T^b$ acts trivially on $H_1(\bar X_n,\Z/a\Z)$, then $b \cdot \tilde x_1  \equiv 0 \mod a$.
Hence  
$b \equiv 0 \mod a$.
\end{proof}

In \cite{Fin13} Proposition~4.4 the result of Proposition~\ref{prop:ParabInXn} is generalised to suitable parabolic elements with positive trace in the Veech group of more general primitive translation surfaces. 

\begin{rem}
 Every parabolic element in $\Gamma(X_n)$ with positive trace is conjugated to a power of $T$. Furthermore, $\Gamma(a)$ is normal for every $a$. 
Thus we know for every parabolic element in $\Gamma(X_n)$ with positive trace, whether it is contained in $\Gamma(a)$ or not. 
\end{rem}

\begin{rem}
\label{rem:torsionfree}
 If we forget about the translation structure on the surface and consider only its topology, then the affine maps can be seen as elements in the mapping class group of the surface. 
 Corollary~1.5 in \cite{I92} tells us that $\Gamma(a) \leq \Gamma(X_n)$ is torsion free for $a \geq 3$. 
 Originally this result is due to Serre (see \cite{Se60}).
 Consequently, $\Gamma(a)$ contains no elliptic elements for $a \geq 3$.
\end{rem}

Remark~\ref{rem:torsionfree} implies that $-I \notin \Gamma(a)$ for $a \geq3$. Thus if $T^k \in \Gamma(a)$ then $-T^k \notin \Gamma(a)$ for $a \geq 3$. 

The knowledge about the parabolic elements in $\Gamma(a)$ is particularly useful if the Veech group of a primitive translation surface is generated by parabolic matrices.
It is well known that this is the case for the surfaces $\bar X_n$ with $n$ odd and $n \geq 5$: 
\begin{lemma}
\label{lemma:XnParabolic}
 For odd $n \geq 5$, 
 $$\Gamma(X_n) = \langle T, R^{-1} T R \rangle\,.$$

 Thus $\Gamma(X_n)$ is generated by parabolic elements. 
\end{lemma}
A proof can be found e.g.\ in \cite{Fin13} Lemma~4.7.

\subsection{Wohlfahrt level}

It is a result due to Wohlfahrt (see \cite{Wo63}) that the \textit{Wohlfahrt level} and the minimal congruence level of every congruence group in $\SL_2(\Z)$ coincide. 
For a finite index subgroup $\Gamma$ in $\SL_2(\Z)$, the Wohlfahrt level is the least common multiple of the cusp widths of the subgroup projected to $\PSL_2(\Z)$. 
A cusp corresponds to a conjugacy class of a maximal parabolic subgroup of $\Gamma$. 
Its width is the index of this subgroup in the ambient maximal parabolic subgroup of $\SL_2(\Z)$.  
Thus the Wohlfahrt level of $\Gamma$ is completely determined by the parabolic elements in $\Gamma$. 

Wohlfahrt's definition relies on the fact that a fundamental domain of $\PSL_2(\Z)$ has exactly one cusp. 
Now consider $\bar X_n$, with $n \geq 5$ odd.
The projection of its Veech group to $\PSL_2(\R)$ is the orientation preserving part of a Hecke triangle group whose fundamental domain in $\H$ has exactly one cusp. 
As $-I_2 \in \Gamma(X_n)$, the projective Veech group is a proper factor group of $\Gamma(X_n)$. 
The parabolic generator $T \in \Gamma(X_n)$ is maximal parabolic, i.e.\ it generates a maximal parabolic subgroup. 
Thus every parabolic element $A \in \Gamma(X_n)$ can be written as $\pm S^{-1} T^m S$ with $m \in \Z$ and $S \in \Gamma(X_n)$. 

In the previous section we learned that $-I_2 \notin \Gamma(a)$ for $a \geq 3$ and that $T^b \in \Gamma(a) \Leftrightarrow a \mid b$, thus it makes sense to use the parabolic elements with positive trace to generalise Wohlfahrt's definition.

\begin{defn}
 Let $A = S^{-1} T^m S$ (with $m \in \Z$ and $S \in \Gamma(X_n)$) be a parabolic element in $\Gamma(X_n)$ with positive trace. 
 We call $|m|$ the \nom{width} of $A$. 
 For a finite index subgroup $\Gamma \leq \Gamma(X_n)$, we define the \nom{(generalised) Wohlfahrt level} of $\Gamma$ as
 \begin{align*}
 \level(\Gamma) \coloneqq \lcm\{ & \wid(A) \mid A = S^{-1} T^m S \in \Gamma \textrm{ where } m \in \Z, S \in \Gamma(X_n) \textrm{ and } \\
 & A \textrm{ is maximal parabolic in } \Gamma \textrm{, i.e.\ } \forall n \mid m : S^{-1} T^n S \notin \Gamma \}.
 \end{align*}
\end{defn}
As $\Gamma$ has finite index in $\Gamma(X_n)$, it contains only finitely many conjugacy classes of maximal parabolic subgroups. 
Thus $\level(\Gamma)$ is well-defined.

In analogy to the work of Wohlfahrt, we define the following subgroups:
\begin{defn}
For $m\geq 1$ let $G(m)$ be the normal closure of $T^ m$ in $\Gamma(X_n)$.
 $$G(m) \coloneqq \nl T^m \nr \subseteq \Gamma(X_n)$$
\end{defn}
We see that the (generalised) Wohlfahrt level of a group $\Gamma$ is densely interwoven with the normal subgroups $G(m)$ contained in $\Gamma$. 

\begin{prop}
 \label{prop:Wohlfahrt1}
 Let $\Gamma \leq \Gamma(X_n)$ be a finite index subgroup with Wohlfahrt level $m$. Then 
 $G(m) \subseteq \Gamma$. If conversely $G(m) \subseteq \Gamma$ then $\level(\Gamma) \mid m$.
\end{prop}
\begin{proof}
 Let $[\Gamma(X_n): \Gamma] = d$. 
 Then for each $A \in \Gamma(X_n)$ there exists a $0 < k \leq d$ such that $A^k \in \Gamma$. 
 We have to show that $S^{-1} T^m S$ belongs to $\Gamma$ for every $S \in \Gamma(X_n)$.
 
 Let $A \coloneqq S^{-1} T S$ and define $l \coloneqq \min \{k > 0 \mid A^k = S^{-1} T^k S \in \Gamma\}$. By the definition of the Wohlfahrt level, $l \mid m$ and therefore $S^{-1} T^m S \in \Gamma$.
 We conclude that $G(m) \subseteq \Gamma$.
 
 Now we want to prove that $G(m) \subseteq \Gamma$ implies that $\level(\Gamma) \mid m$.
 So let $A$ be maximal parabolic with positive trace, then $A=S^{-1} T^k S \in \Gamma$ for some $S \in \Gamma(X_n)$ and $k \in \Z$. We need to show that $k$ divides $m$.
 
 As $G(m) \subseteq \Gamma$, we have $S^{-1} T^m S \in \Gamma$. This implies $S^{-1} T^{\gcd(m,k)} S \in \Gamma$ and as $S^{-1} T^k S$ is maximal parabolic in $\Gamma$, we get that $\gcd(m,k) \geq |k|$, so $\gcd(m,k) = |k|$.
\end{proof}

A simple observation regarding the groups $G(m)$ is the following lemma:
\begin{lemma}
\label{lemma:GmGm'}
 For $m,m' \geq 1$ the product of the groups $G(m)$ and $G(m')$ is
 $$G(m) \cdot G(m') =  \langle G(m), G(m') \rangle = G(\gcd(m, m'))\,.$$ 
\end{lemma}
\begin{proof}
As $\gcd(m, m')$ divides $m$, $G(m) \subseteq G(\gcd(m,m'))$, and analogously $G(m') \subseteq G(\gcd(m,m'))$. 
The groups $G(m)$ and $G(m')$ are normal, thus $ \langle G(m), G(m') \rangle = G(m) \cdot G(m')$. Furthermore, $T^m \in G(m)$ and $T^{m'} \in G(m')$, thus the element $T^{\gcd(m, m')}$ lies in $\langle G(m), G(m') \rangle = G(m) \cdot G(m')$. As $G(m) \cdot G(m')$ is normal, this implies that $G(\gcd(m,m')) \subseteq G(m) \cdot G(m')$.
\end{proof}

\subsection{Comparison of the two level definitions}

Unfortunately the generalised Wohlfahrt level of a congruence group in $\Gamma(X_n)$ does not determine its congruence levels completely (see Example \ref{ex:X5Gamma8}). 
The following theorem makes the correlation between the two level definitions precise:

\begin{thm}
 \label{thm:CongruenceVsWohlfahrtlevel}
 Let $\Gamma \leq \Gamma(X_n)$ be a congruence group, $b$ a minimal congruence level of $\Gamma$ and $a=\level(\Gamma)$ its Wohlfahrt level. 
 Then $a \mid b$ and all prime numbers $p$ dividing $b$ also divide $a$. 
 However, a minimal congruence level of $\Gamma$ does not have to divide the Wohlfahrt level. 
\end{thm}

We prove Theorem~\ref{thm:CongruenceVsWohlfahrtlevel}, with the help of the following lemma. Recall that $\bar{\varphi}_a \colon \Gamma(X_n) \to \Aut((\Z/a\Z)^{2g})$ describes the action of $\Gamma(X_n) \cong \Aff(\bar X_n)$ on $H_1(\bar X_n, \Z/a\Z) \cong (\Z/a\Z)^{2g}$ and that $\Gamma(a) = \ker(\bar \varphi_a)$.

\begin{lemma}
\label{lemma:GGamma=Gamma}
Let $a,b \in \N$ with $\gcd(a,b) = 1$. Then 
\begin{enumerate}[topsep=0cm, parsep=0.1cm, label=\alph*)]
 \item $\bar \varphi_b(G(a)) = \bar \varphi_b(\Gamma(X_n))$ and 
 \item $G(a) \cdot \Gamma(ab) = \Gamma(a)$. 
\end{enumerate}
\end{lemma}
\begin{proof}
a) We have $\Gamma(X_n) \stackrel{\textrm{Lemma } \ref{lemma:XnParabolic}}{=} G(1) = G(\gcd(a,b)) \stackrel{\textrm{Lemma } \ref{lemma:GmGm'}}{=} G(a) \cdot G(b)$.
The normal subgroup $\Gamma(b)$ contains $T^b$ thus $G(b) \subseteq \Gamma(b)$. 
Thus $\bar \varphi_b(\Gamma(X_n)) = \bar \varphi_b(G(a) \cdot G(b)) = \bar \varphi_b(G(a))$, as $G(b) \subseteq \Gamma(b) = \ker(\bar \varphi_b)$. 

b) The subgroups $G(a)$ and $\Gamma(ab)$ are normal in $\Gamma(X_n)$. 
 Thus $G(a) \cdot \Gamma(ab)$ is a group and in particular normal in $\Gamma(X_n)$. 
 As $G(a) \subseteq \Gamma(a)$ and also $\Gamma(ab) \subseteq \Gamma(a)$, we get $G(a) \cdot \Gamma(ab) \subseteq \Gamma(a)$.
 
 For the converse inclusion let $A \in \Gamma(a)$.
 Because of part a) we know that $\bar \varphi_b(A) \in \bar \varphi_b(\Gamma(X_n)) = \bar \varphi_b(G(a))$.
 Thus there exists $B \in G(a)$ such that $\bar \varphi_b(B) = \bar \varphi_b(A)$. Consequently $A$ can be written as
 $A = B \cdot K$,
 where $K \in \Gamma(b) = \ker(\bar \varphi_b)$.
 As $B \in G(a) \subseteq \Gamma(a)$ and $A \in \Gamma(a)$, also $K \in \Gamma(a)$. Hence $K \in \Gamma(a) \cap \Gamma(b)$. But $\Gamma(a) \cap \Gamma(b) = \Gamma(a b)$. This completes the proof.
\end{proof}

Now we can prove that every minimal congruence level of a congruence subgroup in $\Gamma(X_n)$ has only prime divisors that also divide the Wohlfahrt level. This implies in particular that if the Wohlfahrt level is a prime power, then the minimal congruence level is unique. 

\begin{proof}[Proof of Theorem~\ref{thm:CongruenceVsWohlfahrtlevel}]
 Being a congruence group of level $b$ is equivalent to containing $\Gamma(b)$. From $G(b) \subseteq \Gamma(b) \subseteq \Gamma$ and Proposition~\ref{prop:Wohlfahrt1} it follows that $a = \level(\Gamma) \mid b$.
 
 Suppose that $b = c \cdot d$ with $\gcd(a,d)= 1$ and $\gcd(c,d)= 1$. Then by Lemma~\ref{lemma:GGamma=Gamma} 
 $\Gamma(c) = G(c) \cdot \Gamma(cd) = G(c) \cdot \Gamma(b)$.
 As $\gcd(a,b) = \gcd(a,c) \mid c$, it follows that $G(c) \subseteq G(\gcd(a,b))$. 
 Furthermore, $G(a) \subseteq \Gamma$ and $G(b) \subseteq \Gamma(b) \subseteq \Gamma$.
 Thus by Lemma~\ref{lemma:GmGm'} $G(a) \cdot G(b) = G(\gcd(a,b)) \subseteq \Gamma$. 
 Hence $\Gamma(c) = G(c) \cdot \Gamma(b) \subseteq \Gamma$ and $\Gamma$ is a congruence group of level $c$. As $c \mid b$ and $b$ was a minimal congruence level, $d = 1$.
 
 An example which proves that the Wohlfahrt level and the minimal congruence levels can be different follows in Example~\ref{ex:X5Gamma8}.
\end{proof}

\begin{ex}
\label{ex:X5Gamma8}
 The subgroup $U \coloneqq \nl T^4, \Gamma(8) \nr \leq \Gamma(X_5)$ has minimal congruence level $8$ and Wohlfahrt level $4$. 
 
 As $\nl T^4 \nr \subseteq U$, Proposition~\ref{prop:Wohlfahrt1} implies that the Wohlfahrt level of $U$ divides $4 = 2^2$. 
 Theorem~\ref{thm:CongruenceVsWohlfahrtlevel} implies that $U$ has a unique minimal congruence level and that this minimal congruence level is a power of $2$. By construction $\Gamma(8) \subseteq U$, thus the minimal congruence level is at most $8$.
 
 In \cite{Fin13} we give an element $C \in \Gamma(4)$ with $C \notin U$.
 Hence $\Gamma(4) \not\subseteq U$ thus $U$ has minimal congruence level $8$. 
 Alternatively one can verify using a computer algebra system, e.g.\ magma, that $\bar \varphi_8(U) = \bar \varphi_8 (\nl T^4 \nr)$ is a finite group with $32$ elements whereas $\bar \varphi_8(\Gamma(4))$ has $64$ elements. 
 In \cite{Fin13} we also prove that $U \leq \Gamma(X_5)$ is a minimal example with respect to the congruence level for a congruence group in $\Gamma(X_n)$ that has different minimal congruence level and Wohlfahrt level.
\end{ex}

\begin{rem}
 In \cite{Fin13} the Wohlfahrt level is defined more generally for subgroups $\Gamma$ of $\Gamma(X) \subseteq \SL_2(\R)$ for appropriate primitive translation surfaces $\bar X$. 
 The projective Veech group of those surfaces needs to be a Fuchsian group with one cusp. Its Veech group has to contain a parabolic matrix $T$ with positive trace such that $T^b \in \Gamma(a) \Leftrightarrow a \mid b$ and such that $T$ and its conjugates generate $\Gamma(X)$.
 Then Theorem~\ref{thm:CongruenceVsWohlfahrtlevel} can be proven in complete analogy. 
\end{rem}

\subsection{A non-congruence group}

In \cite{Sch13} Gabriela Weitze-Schmithüsen introduces the notion of totally non-congruence groups in $\SL_2(\Z)$. In Theorem~2 she gives a sufficient condition for a finite index subgroup in $\SL_2(\Z)$ to be a totally non-congruence group. 
The next Theorem is a straight forward generalisation of this criterion. Its proof relies in particular on Lemma~\ref{lemma:GGamma=Gamma} a) and thereby on the groups leading to the strong relation between the Wohlfahrt level and the potential congruence levels of a group.

\begin{thm}
\label{thm:NonCongGroup}
 Let $\Gamma \neq \Gamma(X_n)$ be a finite index subgroup in $\Gamma(X_n)$ with Wohlfahrt level $l$. If $l = a \cdot b$ with $\gcd(a,b)=1$ and if $\Gamma$ contains two pairs of parabolic elements with positive trace
 $$(A^{-1} T A)^{m_1}, ((RA)^{-1} T R A)^{m_2} \textrm{ and } (B^{-1} T B)^{m'_1}, ((RB)^{-1} T R B)^{m'_2}$$
 with $A, B \in \Gamma(X_n)$ such that $\lcm(m_1, m_2) \mid a$ and $\lcm(m'_1,m'_2) \mid b$ 
 then 
 $\Gamma$ is not a congruence group in $\Gamma(X_n)$.
\end{thm}
\begin{proof}
 Assume on the contrary that $\Gamma$ is a congruence group with minimal congruence level $c$. Then  we know from Theorem~\ref{thm:CongruenceVsWohlfahrtlevel} that $l \mid c$ and that $c$ can be written as $c = a' \cdot b'$, where $a'$ has the same prime divisors as $a$ and $b'$ has the same prime divisors as $b$. Thus $\gcd(a', b') = 1$.
 
 Our first claim is that the map 
 $$ f \colon \left\{ 
 \begin{array}{rcrcl}
  \bar \varphi_c (\Gamma(X_n)) &\to& \bar \varphi_{a'}(\Gamma(X_n)) &\times& \bar \varphi_{b'}(\Gamma(X_n)) \vspace{0.1cm}\\ 
  \bar \varphi_c(A) &\mapsto& (\;\bar \varphi_{a'}(A)&,& \bar \varphi_{b'}(A)\;)
 \end{array}
 \right.$$ 
 is a group isomorphism. 
 As $\ker(\bar\varphi_c) = \Gamma(c) \subseteq \Gamma(a') = \ker(\bar \varphi_{a'})$ and $\Gamma(c) \subseteq \Gamma(b')$, the map is a well defined group homomorphism.
 It is injective, because it is the restriction of the canonical group isomorphism 
 $\SL_{n-1}(\Z/c\Z) \cong \SL_{n-1}(\Z/a'\Z) \times \SL_{n-1}(\Z/b'\Z)$ to $\bar \varphi_c (\Gamma(X_n))$.
 Thus it remains to prove the $f$ is surjective.
 So let $A \in \bar\varphi_{a'}(\Gamma(X_n))$ and $B \in \bar\varphi_{b'}(\Gamma(X_n))$. 
 As $\gcd(a', b') = 1$ Lemma~\ref{lemma:GGamma=Gamma} a) implies $\bar \varphi_{a'}(G(b')) = \bar \varphi_{a'}(\Gamma(X_n))$ and $\bar \varphi_{b'}(G(a')) = \bar \varphi_{b'}(\Gamma(X_n))$. 
 Thus there exists $\tilde A \in G(b')$ and $\tilde{B} \in G(a')$ such that $\bar\varphi_{a'}(\tilde A) = A$ and $\bar\varphi_{b'}(\tilde B) = B$.
 We know that for all $q \geq 2: G(q) \subseteq \Gamma(q)$.
 It follows that $\bar\varphi_{a'}(\tilde{A} \tilde{B}) = \bar\varphi_{a'}(\tilde{A}) = A$ and $\bar\varphi_{b'}(\tilde{A} \tilde{B}) = \bar\varphi_{b'}(\tilde{B}) = B$. 
 We deduce $f(\bar\varphi_{c}(\tilde{A} \tilde{B})) = (A,B)$.
 
 As $\gcd(a', b') = 1$, there exist $l,k \in \Z$ such that $l a' + k b' = 1$. 
 Because of $\lcm(m_1, m_2) \mid a$ and $a \mid a'$ there exist $h_1,h_2 \in \Z$ such that 
 $m_1 h_1 = m_2 h_2 = a'$.
 We know that $(A^{-1} T A)^{m_1}, ((RA)^{-1} T R A)^{m_2} \in \Gamma$ thus 
 $(A^{-1} T A)^{m_1 h_1 l} \in \Gamma$ and $((RA)^{-1} T R A)^{m_2 h_2 l} \in \Gamma$.
 It follows that 
 \begin{align*}
  \bar\varphi_{c}((A^{-1} T A)^{m_1 h_1 l})
  & = (\, \bar\varphi_{a'}(\,(A^{-1} T A)^{la'}\,) \,,\, \bar\varphi_{b'}(\,(A^{-1} T A)^{la'} \,) \,)\\  
  & = (\, I \,,\, \bar\varphi_{b'}(\,(A^{-1} T A)^{la'} \cdot (A^{-1} T A)^{kb'} \,) \,) \\
  & = (\, I \,,\, \bar\varphi_{b'}(A^{-1} T A) \,) \in \bar\varphi_{c}(\Gamma)
 \end{align*}
 and analogously 
 $\bar\varphi_{c}(((RA)^{-1} T RA)^{m_2 h_2 l})
  = (\, I \,,\, \bar\varphi_{b'}((RA)^{-1} T RA) \,) \in \bar\varphi_{c}(\Gamma)$.
Lemma~\ref{lemma:XnParabolic} states that $T$ and $R^{-1} T R$ generate $\Gamma(X_n)$. Thus also $A^{-1} T A$ and $(RA)^{-1} T R A$ generate $\Gamma(X_n)$. This implies that 
$$\{I\} \times \bar\varphi_{b'}(\Gamma(X_n)) \subseteq \bar\varphi_{c}(\Gamma) \,.$$ 
Using the same arguments on $(B^{-1} T B)^{m'_1 h'_1 k}$ and $((RB)^{-1} T R B)^{m'_2 h'_2 k}$ where $h'_1, h'_2$ are chosen such that $m'_1 h'_1 = m'_2 h'_2 = b'$ we obtain 
$$\bar \varphi_{a'}(\Gamma(X_n)) \times \{I\} \subseteq \bar \varphi_{c}(\Gamma)\,.$$
Altogether this proves  
$\bar \varphi_c(\Gamma(X_n)) = \bar \varphi_{a'}(\Gamma(X_n)) \times \bar\varphi_{b'}(\Gamma(X_n)) = \bar \varphi_{c}(\Gamma)$. 
But as $c$ was supposed to be a congruence level of $\Gamma$, this implies $\Gamma = \Gamma(X_n)$ and thereby contradicts our assumption. 
\end{proof}

In the end we give an example for a subgroup $\Gamma \leq \Gamma(X_5)$ that fulfils the requirements of Theorem~\ref{thm:NonCongGroup} and is thereby an example of a non-congruence group.

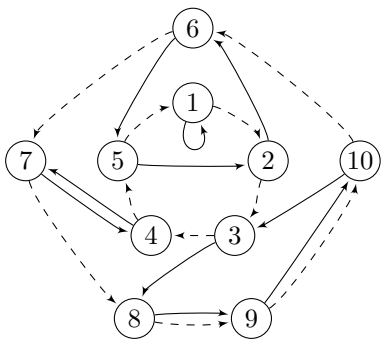
\begin{figure}[htbp]
 \centering
 \begin{tikzpicture}[scale=1.1]

\tikzstyle{place}=[draw,circle,inner sep=1pt,minimum size=15pt];
\tikzstyle{basepoint}=[draw,rectangle,inner sep=2pt];
\tikzstyle{myarrow}= [->, >=latex', shorten >=1pt]; 

\node[place] (1) at (0,0) [] {$1$};
\node[place] (2) at (0.9,-0.7) [] {$2$};
\node[place] (3) at (0.5,-1.6) [] {$3$};
\node[place] (4) at (-0.5,-1.6) [] {$4$};
\node[place] (5) at (-0.9,-0.7) [] {$5$};

\node[place] (6) at (0,0.9) [] {$6$};
\node[place] (7) at (-2,-0.7) [] {$7$};
\node[place] (8) at (-0.7,-2.6) [] {$8$};
\node[place] (9) at (0.7,-2.6) [] {$9$};
\node[place] (10) at (2,-0.7) [] {$10$};

\draw [myarrow] (1) ..controls +(250:0.7cm) and +(290:0.7cm)  .. (1);

\draw [myarrow] (2) ..controls +(90:0.4cm) and +(-30:0.4cm)  .. (6);
\draw [myarrow] (6) ..controls +(210:0.4cm) and +(90:0.4cm)  .. (5);
\draw [myarrow] (5) ..controls +(-10:0.4cm) and +(190:0.4cm)  .. (2);

\draw [myarrow] (4) ..controls +(150:0.4cm) and +(-20:0.4cm)  .. (7);
\draw [myarrow] (7) ..controls +(-40:0.4cm) and +(170:0.4cm)  .. (4);

\draw [myarrow] (3) ..controls +(200:0.4cm) and +(70:0.4cm)  .. (8);
\draw [myarrow] (8) ..controls +(10:0.4cm) and +(170:0.4cm)  .. (9);
\draw [myarrow] (9) ..controls +(50:0.4cm) and +(240:0.4cm)  .. (10);
\draw [myarrow] (10) ..controls +(220:0.4cm) and +(20:0.4cm)  .. (3);

\draw [myarrow,dashed] (1) ..controls +(-10:0.4cm) and +(110:0.4cm)  .. (2);
\draw [myarrow,dashed] (2) ..controls +(250:0.4cm) and +(50:0.4cm)  .. (3);
\draw [myarrow,dashed] (3) ..controls +(180:0.4cm) and +(0:0.4cm)  .. (4);
\draw [myarrow,dashed] (4) ..controls +(130:0.4cm) and +(290:0.4cm)  .. (5);
\draw [myarrow,dashed] (5) ..controls +(70:0.4cm) and +(190:0.4cm)  .. (1);

\draw [myarrow,dashed] (10) ..controls +(110:0.5cm) and  +(-10:0.5cm) .. (6);
\draw [myarrow,dashed] (9) ..controls +(30:0.5cm) and +(260:0.5cm) .. (10);
\draw [myarrow,dashed] (8) ..controls +(-10:0.5cm) and +(190:0.5cm) .. (9);
\draw [myarrow,dashed] (7) ..controls +(280:0.5cm) and +(140:0.5cm) .. (8);
\draw [myarrow,dashed] (6) ..controls +(190:0.5cm) and +(70:0.5cm) .. (7);

\end{tikzpicture}
 \caption{Coset graph of a non-congruence subgroup $\Gamma \leq \Gamma(X_n)$.}
 \label{fig:NonCongEx}
\end{figure}

\begin{ex}
Figure \ref{fig:NonCongEx} shows the left coset graph of the subgroup 
$$ \Gamma = \langle T, R^{-1}TR^{-1}, R T^2R, R^2TRTR, R^{-2}T^{-2}RTR^2, R^{-2}T^2RTR^2 \rangle \leq \Gamma(X_5)\,,$$
i.e.\ a graph whose vertices are the left cosets of $\Gamma(X_5)/\Gamma$, 
solid edges $A \cdot \Gamma \to T A \cdot \Gamma$ and dashed edges $A \cdot \Gamma \to R A \cdot \Gamma$. The coset $I \cdot \Gamma$ carries the label $1$.

A presentation of $\Gamma(X_5)$ is $\langle R,T \mid R^{10} = I, (T^{-1}R)^2 = R^5, R^5 T = T R^5 \rangle$. 
It is easy to check that $R^5$ and $(T^{-1}R)^2$ describe a closed path at each vertex of the coset graph. Consequently the defining relations of $\Gamma(X_5)$ are fulfilled and Figure~\ref{fig:NonCongEx} shows indeed the coset graph of a subgroup of $\Gamma(X_5)$. 
From the coset graph one easily reads off that $T,\, R^{-1}T^3R \in \Gamma$ as well as $R^{-2} T^4 R^2,\, R^{-3} T^2 R^3 \in \Gamma$. 
Hence we can apply Theorem~\ref{thm:NonCongGroup} with $A=I$, $B=R^2$, $m_1=1$, $m_2=3$, $m'_1=4$ and $m'_2=2$. 
We conclude that $\Gamma$ is an example for a subgroup of $\Gamma(X_n)$ that is not a congruence group.
\end{ex}

\printbibliography

\end{document}